\documentclass{article}
\usepackage[utf8]{inputenc}
\usepackage[dvips]{graphics,graphicx}
\usepackage{amsfonts,amssymb,amsmath,color,mathrsfs, amstext}
\usepackage{amsbsy, amsopn, amscd, amsxtra, amsthm,authblk,enumerate,algorithmicx,algorithm}
\usepackage{bm,bbm}
\usepackage{algpseudocode}
\usepackage{upref}
\usepackage{geometry}
\usepackage{ulem}
\usepackage{float}
\usepackage{yhmath}
\usepackage[colorlinks,
            linkcolor=red,
            anchorcolor=red,
            citecolor=red
            ]{hyperref}

\numberwithin{equation}{section}

\newcommand{\abs}[1]{\left\vert#1\right\vert}

\def\cO{\mathcal{O}}

\newtheorem{theorem}{Theorem}[section]

\newtheorem{lemma}{Lemma}[section]
\newtheorem{proposition}{Proposition}[section]

\numberwithin{equation}{section}
\numberwithin{figure}{section}
\numberwithin{table}{section}
\DeclareMathOperator{\sgn}{sgn}
\def\half{\frac 1 2}

\begin{document}

\title{Some random batch particle methods for the Poisson-Nernst-Planck and Poisson-Boltzmann equations}

\author[1]{Lei Li\thanks{leili2010@sjtu.edu.cn}}
\author[2]{Jian-Guo Liu\thanks{jliu@phy.duke.edu}}
\author[3]{Yijia Tang\thanks{yijia\underline{~}tang@sjtu.edu.cn}}
\affil[1]{School of Mathematical Sciences, Institute of Natural Sciences, MOE-LSC, Shanghai Jiao Tong University, Shanghai, 200240, P. R. China.}
\affil[2]{Department of Mathematics and Department of Physics, Duke University, Durham, NC 27708, USA.}
\affil[3]{School of Mathematical Sciences, Shanghai Jiao Tong University, Shanghai, 200240, P. R. China.}

\date{}
\maketitle
\begin{abstract}

We consider in this paper random batch interacting particle methods for solving the Poisson-Nernst-Planck (PNP) equations, and thus the Poisson-Boltzmann (PB) equation as the equilibrium, in the external unbounded domain. To justify the simulation in a truncated domain, an error estimate of the truncation is proved in the symmetric cases for the PB equation. Then, the random batch interacting particle methods are introduced which are $\cO(N)$ per time step. 
The particle methods can not only be considered as a numerical method for solving the PNP and PB equations, but also can be used as a direct simulation approach for the dynamics of the charged particles in solution. 
The particle methods are preferable due to their simplicity and adaptivity to complicated geometry, and may be interesting in describing the dynamics of the physical process.
Moreover, it is feasible to incorporate more physical effects and interactions in the particle methods and to describe phenomena beyond the scope of the mean-field equations. 
\end{abstract}


\section{Introduction}
\label{sec: intro}

The charge distribution in dilute ionic solution around some charged surfaces is important for a wide range of applications in electrochemistry \cite{Gouy,Chapman}, biophysics \cite{HonigClassical,Davis1990Electrostatics} and colloidal physics \cite{DeLa41,VeOv48}. In the so-called implicit solvent model, the solvent is modeled by a continuum while the ions can be either modeled by charged particles or continuum distributions. 
When one models the ions using continuum distributions, some partial differential equations (PDEs) can be proposed. 
Regarding the ion transport, the Poisson-Nernst-Planck (PNP) equations \cite{EiHyLiu10,Flavell2014A,Ji_et_all18,LiuEisenberg2020} have been used to describe the non-equilibrium processes in the dilute regime.
The Poisson-Boltzmann (PB) equation, proposed by Gouy \cite{Gouy} and Chapman \cite{Chapman} independently, can be viewed as the equilibrium of the PNP equations. The PB equation is a typical implicit solvent model to describe the distribution of the electric potential in dilute solution at equilibrium state when an object with free charges inside is immersed into an ionic solution.
Various numerical methods have been proposed for the PNP and PB equations based on the PDE descriptions in literature \cite{Lu2008Recent,cai2013}, such as the finite difference method (FD) \cite{JCC1988,Chern2010Accurate,Flavell2014A}, finite element method (FEM) \cite{Baker_et_al2000,ChenHolstXu07} and boundary element/integral method \cite{Bordner2003Boundary,Lu19314}.

While the continuum descriptions of the charge distributions can capture some mean-field behaviors, the numerical simulations based on particles, or molecular dynamics (MD) simulations, can potentially give more physics and give some dynamical properties of the systems \cite{MD}. The MD simulations with electrostatic Coulomb interactions are usually challenging due to the long-range nature. A lot of efforts have been made to efficiently approximate the pair-wise interactions between charges in an electrolyte, such as the fast multipole method (FMM) \cite{FMM}, the Ewald method \cite{Ewald}, particle mesh Ewald (PME) \cite{PME} and particle-particle particle mesh Ewald (PPPM) \cite{PPPM}. Recently, a stochastic method, the Random batch Ewald (RBE) method \cite{RBE}, was proposed to simplify particle simulations with Coulomb interaction. 
FMM can reduce the cost to $\mathcal{O}(N)$ per time step but the implementation is nontrivial and the efficiency can be observed when the number of particles is large. 
The Ewald-based methods like PPPM or RBE methods can reduce the cost to $\cO(N\log N)$ or $\cO(N)$ per time step, but the simulations take place in a box with periodic boundary conditions (BCs). Another popular method for plasma simulations is the particle-in-cell (PIC) method \cite{PIC}. The PIC method considers the interaction between particles by computing the electric field on a deterministic grid and coupling the charged particles to the field, which has a cost of $\cO(N\log N)$ using FFT. The collision-field method \cite{1998Improved} can be viewed as an improved version of PIC. It treats the inter-species collisions in deterministic grid like PIC, while treats the intra-species collisions through the Langevin equation. This method can ensure momentum and energy conservation exactly using velocity corrections.

In this work, we would like to seek some particle methods for the PNP and PB equations using the random batch idea. 
In our particle methods, we simulate the overdamped Langevin equations which are the microscopic descriptions corresponding to the PNP equations so that the distributions in large time will solve the PB equation. In the simulation, each particle will interact with the others through the long-range Coulomb interaction. Note that solving the PDEs using particle methods with interaction is often expensive and the accuracy is lower compared with solving the PDEs directly.  We emphasize the advantages of using particle methods in several aspects. On the one hand, particle simulation can render the transient phenomena and some features that the mean-field PDEs cannot capture, though in this work our main focus is to solve the PDEs and the equilibrium distributions. On the other hand, the particle method is meshfree so that is insensitive to dimensionality and geometry. Moreover, thanks to the random batch method (RBM) for interacting particle systems proposed by Jin et al in \cite{JinLiLiu20}, we can reduce the simulation cost per iteration to $\cO(N)$ so that the computational cost can be addressed to some extent. The RBM utilizes small but random batch idea so that interactions only occur inside the small batches to reduce the computational cost per time step from $\mathcal{O}(N^2)$ to $\mathcal{O}(N)$ in a surprisingly simple way.
	
In the original setting, the charged object is immersed in some unbounded solution. For the particle simulation, one must truncate the domain and prescribe a suitable BC to the artificial surface. We choose to use the reflecting BC and will provide an error estimate (Theorem \ref{thm: convergence}) for the truncation in the PDE level which says the solution of the truncated PB equation can well approximate the solution of the original PB equation. For symmetric cases, the $L_1$ error of the solution in the truncated domain with an extra Neumann BC and the solution in the unbounded domain is exponentially decaying in the truncated length. This means one only needs to simulate the ionic particles in a suitably truncated domain. 

In the mean-field PDE descriptions of the ionic solutions (i.e., the PNP and PB equations), only the effects of Coulomb interactions are present and the hard sphere potentials like the Lennard-Jones potential \cite{LJ1924} are ignored as the sizes of the particles tend to zero in the mean-field limit. 
A naive particle method for the PNP and PB equations will thus include the Coulomb forces only for the simulation. Such systems are often troublesome: the positive charges and negative charges can merge if there is solely Coulomb interaction. This brings difficulty to direct particle simulations. To address this issue, one often includes the physical hard sphere potentials so that the simulation will be meaningful. 
We remark, however, that if we solely want to capture the mean-field behaviors described by the PNP and PB equations, the RBM can resolve the issue of attraction between opposite charges so that we do not need the hard sphere potentials.  The reason is that if we do random reshuffling of particles at each time step, there is no possibility that two particles stay in the same batch all the time. Hence, as the number of particles tends to infinity, the simulation results can capture the mean-field behaviors. This means applying RBM to the $N$-particle system with only Coulomb interaction can serve as a {\it numerical particle method} for the PNP and PB equations.
Of course, if one wants to do MD simulations where the hard sphere potentials have physical impacts, one must consider the hard sphere potentials and it is indeed relatively easy to include these effects in our particle method. Our particle methods can thus also be used as MD simulation approaches with these effects included.  
For instance, we conduct a colloidal example in section \ref{subsec: colloid}, where the charge reversal phenomenon can be observed in accordance with experiments \cite{experiment}, theories \cite{RevModPhys} and MD simulations \cite{colloid,RBE}, which can not be described by the mean-field PB equation.

Let us give here some comments on the random batch particle methods for the comparison to the aforementioned methods.
PIC is traditionally a popular method for PB equation and plasma simulations. Compared to PIC, our particle methods discretize the overdamped Langevin equations, and consider interactions between particles directly. Therefore, various physical features can be added in the simulation, including the finite size effect and additional interactions such as Lennard-Jones \cite{LJ1924}.
The random batch particle methods are based on Monte Carlo ideas so the accuracy may not be very high when the time step is not small or the batch size is not big compared to PIC or FMM. However, they are mesh-free so that they are easy to implement, and the prefactor in the linear scaling is smaller compared to FMM. Moreover, the random batch particle methods can be better suited for parallel computation. In fact, the super-scalability based on the random batch idea has already been confirmed in the RBE for MD \cite{liang2021superscalability}.

The rest of the paper is organized as follows. In section \ref{sec: setup}, we briefly introduce the PNP and PB equations from a microscopic point of view and some basic setup. In section \ref{sec: model}, we present an approximate model by truncating the domain and we obtain an error estimate in section \ref{sec: approx}, which demonstrates that the solution of the approximate PB equation converges exponentially to the solution of the original PB equation in some special cases. In section \ref{sec: RBM}, we explain the details of the random batch particle methods for the PNP and PB equations in a truncated domain and address several important issues in implementation. Numerical examples are given in section \ref{sec: numerics}, which show the effectiveness of the random batch particle methods in solving the PNP and PB equations.
Concluding remarks are drawn in section \ref{sec: Conclusions}.

\section{The PNP and PB equations}
\label{sec: setup}

In this section, we first give a brief derivation of the PNP and PB equations starting from the Langevin equation describing the motion of the microscopic particle. Then we perform nondimensionalization for them.

\subsection{The mathematical setup}

\begin{figure}[H]
\centering
\includegraphics[width=0.5\textwidth]{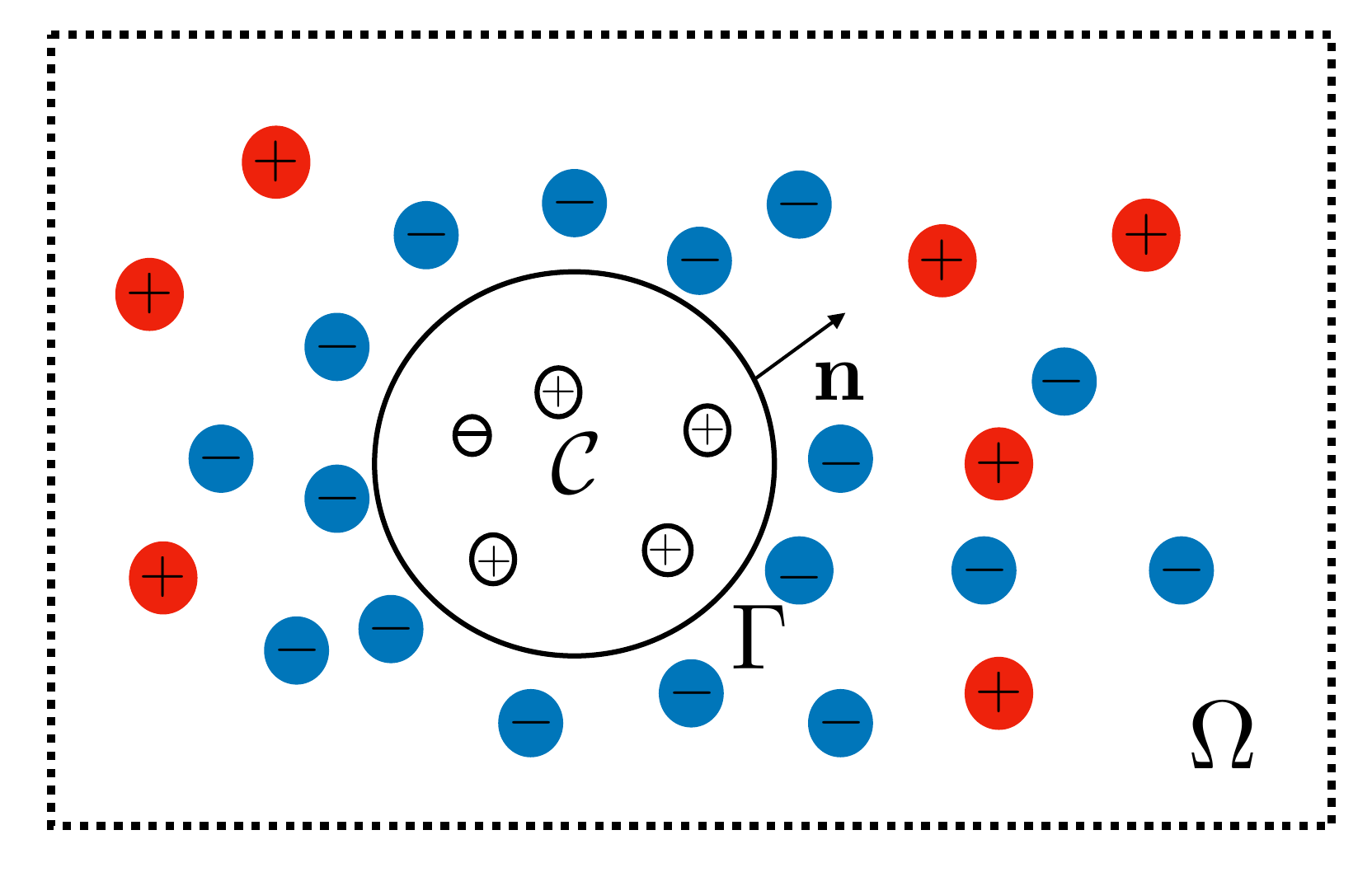}
\caption{Illustration of the domain.}
\label{fig: illustr}
\end{figure}

Consider an object $\mathcal{C}\subset \mathbb{R}^d$ with some free charges inside, immersed in some electrolyte solution $\Omega$ with $J$ ionic species. Here, the object $\mathcal{C}$ could be a macromolecule or the solute in solvent. It could also be a cell (like a neuron) in the tissue fluid.  In this setup, $\Omega=\mathbb{R}^d\backslash \bar{\mathcal{C}}$.  Let $\rho_f$ be the distribution of charges inside $\mathcal{C}$ which we regard as given in our consideration. Denote $\Gamma=\partial \mathcal{C}$, the boundary of $\mathcal{C}$. This is either the macromolecule/solvent-solute interface or the membrane of the cell. Later, we will generally call this the ``interface". We will also assume that the interface $\Gamma$ has a fixed shape and the ions can not go through it. Then the ions will concentrate close to the interface and form a screening layer in several Debye lengths. The setup is illustrated in Figure \ref{fig: illustr}.
We assume the dielectric constant
does not change too much from inside to outside of $\mathcal{C}$, so we will assume the dielectric constant $\varepsilon$ is the same in $\mathbb{R}^d$ throughout this work to illustrate our particle methods. Note that this assumption is reasonable for the cell immersed in tissue fluid but quite restrictive for the solvent-solute models. However, for the solvent-solute models, the apparent surface charge is often known like the colloid immersed in a solution. Hence, one can equivalently use an effective $\rho_f$ with the same dielectric constant inside $\mathcal{C}$. In some applications, the effects of a different dielectric constant inside $\mathcal{C}$ should be considered in detail. For these cases, we need to make use of the image charge method to apply our particle methods \cite{Cai09,PhysRevE.87.013307}. See Appendix \ref{sec: dielectric} for some discussions on the effects of different dielectric constants. 

In this work, the finite size of the ions outside $\mathcal{C}$ is considered negligible, so the ions can be treated as point charges. For a typical ionic particle $X_j$ of the $j$-th species $(1\leq j\leq J)$ with valence $z_j$, it is subject to the electric field and collides with other particles and molecules of the solvent. We assume the electrolyte is dilute so that the collision between particles and collision with the solvent molecules are modeled by friction and white noise. 
Then the motion of $X_j$ is described by the overdamped Langevin system
\begin{equation}
\label{eq: overdamped Langevin}
	d X_j=-\frac{1}{\gamma}\nabla(z_je\Phi)dt+\sqrt{2D}\, dB+dR.
\end{equation}
Here $dX_j$ comes from the friction term so that  $X_j\in \mathbb{R}^d$ is the location of the particle. $e$ is the electron charge, $\Phi$ is the total electronic potential, $-z_je\nabla \Phi$ represents the electronic force. $D$ is the diffusion constant satisfying Einstein relation $D=k_B T/\gamma$, where $k_B$, $T$ and $\gamma$ denote the Boltzmann constant, the absolute temperature and the viscosity coefficient respectively. 
Note that here $[B]=\sqrt{[t]}$, where $[\cdot]$ represents the dimension, (the scale of time $[t]$ is chosen as $\dfrac{L_c^2}{D}$ in the next subsection), this gives $B(t)=\sqrt{[t]}\tilde{B}\left(\frac{t}{[t]}\right)$ with $\tilde{B}(\cdot)$ being a standard Brownian motion. Furthermore, $R(t)\in \mathbb{R}^d$ is the reflecting process associated to $X_j$ which prohibits $X_j$ from crossing the interface. Then $X_j(t)$ is a $\bar{\Omega}$-valued process and $R(t)$ satisfies
 \begin{equation}
\label{eq: reflection process}
	R(t)=\int_0^t\bm{n}(X_j)\mathrm{d}|R|_s,\quad |R|_t=\int_0^t \mathbbm{1}_{\Gamma}(X_j)\mathrm{d}|R|_s.
\end{equation}
Here, $\bm{n}(X_j)$ denotes the outward unit normal to $\Gamma$ at the point $X_j\in \Gamma$, $R(0)=0$, and $|R|_t$ is the total variation of $R(t)$ on $[0,t]$, namely,
$|R|_t=\sup \sum_{k=1}^n\left|R(t_k)-R(t_{k-1})\right|,$
where the supremum is taken over all partitions $0\leq t_0<t_1<\cdots<t_n=t$. 

The reflecting stochastic differential equation (RSDE)  \eqref{eq: overdamped Langevin}-\eqref{eq: reflection process} is also called Skorokhod SDE. Finding the solution pair $(X_j, R)$ to the RSDE is the well-known Skorokhod problem, which was pioneered by Skorokhod in \cite{Skorokhod1961}, where the reflecting process in a half-line $[0,\infty)$ was taken into account. Later on, Tanaka \cite{Tanaka79} considered the multi-dimensional case in a convex domain. Lions and Sznitman \cite{LionsSznitman84} considered more a general domain satisfying the admissibility condition.
At the discrete level, RSDE can be solved by some standard numerical methods \cite{LEPINGLE1995119,Pettersson95}, which combine the widely used Euler-Maruyama scheme and some reflecting techniques.

Next, we turn to the continuous description. Let $\rho_j$ be the macroscopic species density of the $j$-th species in $\Omega$ (outside the $\mathcal{C}$). By It\^o's formula, $\rho_j$ corresponding to systems \eqref{eq: overdamped Langevin}-\eqref{eq: reflection process} is governed by the Nernst-Planck equation 
\begin{equation}
\label{eq: NP}
\left\{
\begin{array}{ll}
	\partial_t \rho_j= \nabla\cdot\left[D\left(\nabla\rho_j+\frac{z_je}{k_BT}\rho_j\nabla\Phi\right)\right],\quad \bm{x}\in\Omega, t>0,\\[3mm]
    \rho_j|_{t=0}=\rho_{0,j},\quad \bm{x}\in\Omega,\\[3mm]
\left\langle D\left(\nabla\rho_j+\frac{z_je}{k_BT}\rho_j\nabla\Phi\right),\bm{n}\right\rangle=0,\quad \bm{x}\in \Gamma,
\end{array}
	\right.
\end{equation}
coupled with the Poisson equation for the total electronic potential $\Phi$
\begin{equation}
\left\{
\begin{array}{lll}
    -\varepsilon \Delta \Phi=e\rho_f ,\quad \bm{x}\in \mathcal{C},\\[2mm]
	-\varepsilon \Delta \Phi=\sum\limits_{j=1}^Jz_je\rho_j,\quad \bm{x}\in \Omega,\\[3.5mm]
	[\Phi]\big|_{\Gamma}=0, \quad 
	[\nabla \Phi\cdot \bm{n}]\big|_{\Gamma}=0.
\end{array}\right.
  \label{eq:Poisson}
\end{equation}
\eqref{eq:Poisson} is derived from Gauss's law. Here, $\varepsilon=\varepsilon_0\varepsilon_r$ is the dielectric constant
both inside and outside the $\mathcal{C}$, $\bm{n}$ is the outward unit normal of $\Gamma$ (pointing to $\Omega$).
The coupling of the Nernst-Planck equation \eqref{eq: NP} and the Poisson equation \eqref{eq:Poisson}  is the well-known PNP system.

For the Nernst-Planck equation \eqref{eq: NP}, the stationary distribution has the form of
$$
\rho_j=c_j\exp\left(-\dfrac{z_je}{k_BT}\Phi\right)
$$ 
with $c_j$ a positive constant.
Provided the system satisfies the electroneutrality condition
$$\sum_{j=1}^J z_j\rho_j^{\infty}=0,$$
where $\rho_j^{\infty}$ is the far field concentration, one has $\Phi(\bm{x})\rightarrow 0$ as $ |\bm{x}| \rightarrow \infty$ and thereby $c_j=\rho_j^{\infty}$.  
Hence, the stationary species density is given by the Boltzmann distribution
\begin{equation}
	\rho_j=\rho_j^{\infty}\exp\left(-\dfrac{z_je}{k_BT}\Phi\right),\quad j=1,\cdots,J.
\end{equation}
Consequently, the steady state Poisson equation becomes
\begin{equation}
\label{eq: PB-BVP}
\left\{
\begin{array}{lll}
	-\varepsilon \Delta \Phi=e\rho_f ,\quad \bm{x}\in \mathcal{C},\\[1.5mm]
	-\varepsilon\Delta \Phi=\sum\limits_{j=1}^Jz_je\rho_j^{\infty}\exp\left(-\dfrac{z_je}{k_BT}\Phi\right),\quad \bm{x}\in \Omega,\\[3.5mm]
[\Phi]\big|_{\Gamma}=0, \quad [\nabla \Phi\cdot \bm{n}]\big|_{\Gamma}=0,
	\quad \Phi(\bm{x})\rightarrow 0 \ \text{as} \ |\bm{x}| \rightarrow \infty.
\end{array}
	\right.
\end{equation}
This is the so-called PB equation. 

\subsection{Nondimensionalization}

Denote $L_c$ as the diameter of $\mathcal{C}$ and $\rho_c$ as the characteristic concentration. Introduce the Debye length defined by 
$\lambda_D=\sqrt{\dfrac{\varepsilon k_B T}{e^2\rho_c}}$ and a parameter $\nu=\left(\dfrac{\lambda_D}{L_c}\right)^2$.
Then, one can rescale the variables 
$\tilde{\bm{x}}=\dfrac{\bm{x}}{L_c}, \tilde{t}=\dfrac{Dt}{L_c^2}$, 
let $\tilde{\mathcal{C}}=\{\tilde{\bm{x}}\in \mathbb{R}^d: \tilde{\bm{x}}L_c\in \mathcal{C}\},\tilde{\Gamma}=\partial \tilde{\mathcal{C}}, \tilde{\Omega}=\mathbb{R}^d\backslash \overline{\tilde{\mathcal{C}}}$
and introduce the following dimensionless quantities 
$$
\tilde{\Phi}=\dfrac{e\Phi}{k_BT},
\quad \tilde{\rho}_j=\frac{\rho_j}{\rho_c},\quad \tilde{\rho}_j^{\infty}=\frac{\rho_j^{\infty}}{\rho_c}, \quad \tilde{\rho}_f=\frac{\rho_f}{\rho_c},\quad \tilde{\rho}_{0,j}=\frac{\rho_{0,j}}{\rho_c}.
$$
For notation convenience, the tildes over all quantities are dropped from now on. 

For simplicity, we only consider the symmetric monovalent electrolyte (we also consider the asymmetric case in a numerical example) throughout this paper, i.e. $J=2, j=\pm, z_+=1, z_-=-1$.
The dimensionless PNP system for two species reads
\begin{equation}
\label{eq:nondim PNP}
\left\{
\begin{array}{llllll}
	\partial_t \rho_j =
	\nabla\cdot\left(\nabla\rho_j+z_j\rho_j\nabla \Phi\right),\quad \bm{x} \in  \Omega,\\[2mm]
	\langle \nabla\rho_j+z_j\rho_j\nabla \Phi,\bm{n}\rangle=0,\quad \bm{x} \in  \Gamma,\quad j=\pm,\\[2.5mm]
	\rho_j|_{t=0}=\rho_{0,j},\quad \bm{x}\in\Omega,\\[2mm]
	-\nu \Delta \Phi=\rho_f ,\quad \bm{x}\in \mathcal{C},\\[1.5mm]
	-\nu\Delta \Phi=\rho_{+}-\rho_{-},\quad \bm{x} \in \Omega,\\[2.5mm]
	[\Phi]\big|_{\Gamma}=0, \quad [\nabla \Phi\cdot \bm{n}]\big|_{\Gamma}=0,
	\quad \Phi(\bm{x})\rightarrow 0 \ \text{as} \ |\bm{x}| \rightarrow \infty,
\end{array}
	\right.
\end{equation}
while the equilibrium distributions read
$$
\rho_{j}=\rho_{j}^{\infty}\exp\left(-z_j\Phi\right),\quad j=\pm.$$
Due to electroneutrality, 
$$\rho_{\infty}:=\rho_{+}^{\infty}=\rho_{-}^{\infty}.$$
Then, the nonlinear dimensionless PB system under investigation is
 \begin{eqnarray}
 \label{eq:PB}
\left\{
\begin{array}{lll}
-\nu \Delta \Phi=\rho_f ,\quad \bm{x}\in \mathcal{C},\\[1.5mm]
-\nu\Delta \Phi=\rho_{\infty}\left(e^{-\Phi}-e^{\Phi}\right),\quad \bm{x} \in \Omega,\\[2.5mm]
	[\Phi]\big|_{\Gamma}=0, \quad [\nabla \Phi\cdot \bm{n}]\big|_{\Gamma}=0,
	\quad \Phi(\bm{x})\rightarrow 0 \ \text{as} \ |\bm{x}| \rightarrow \infty.
\end{array}
	\right.
\end{eqnarray}
From formal asymptotic expansion, we know there is a boundary layer (the so-called Debye screening layer) around $\Gamma$ in the solvent region $\Omega$ whose thickness is of $\mathcal{O}(\sqrt{\nu})$ \cite{Allaire_2013}.

\section{Approximation by the truncated domain}
\label{sec: model}
As we mentioned, we seek particle methods that directly simulate the dynamics of the ions. At the microscopic level, the underlying RSDE for \eqref{eq:nondim PNP} is given by
\begin{equation}
\label{eq: nondim Langevin}
\left\{
\begin{array}{ll}
	dX_j=-\nabla\left(z_j\Phi\right)dt+\sqrt{2}dB+dR,\\[3mm]
	R(t)=\int_0^t\bm{n}(X_j)\mathrm{d}|R|_s,\quad |R|_t=\int_0^t \mathbbm{1}_{\Gamma}(X_j)\mathrm{d}|R|_s.
	\end{array}
\right.
\end{equation}
This is the dimensionless version of \eqref{eq: overdamped Langevin}. Again, the reflecting process $R$ associated with $X_j$ prevents $X_j$ from going into $\mathcal{C}$.

However, it is unrealistic to do the simulation in an unbounded domain. Due to no-flux BC on $\Gamma$, we have $\int_{\Omega}\rho_{\pm}\mathrm{d}x =\int_{\Omega}\rho_{\infty}e^{\mp\Phi}\mathrm{d}x=+\infty$ as $\rho_{\infty}\neq 0$, which means the total positive and negative charges are infinite in the unbounded external domain. However, we can only simulate finite number of particles in particle simulation. In this section, we propose a model in the truncated domain for approximation.

\subsection{RSDE with an artificial wall} 
\label{subsec: RSDE in bounded domain}

Intuitively, consider, for example, $\mathrm{NaCl}$ solution in an unbounded container. Put a charged object into it, a screening layer is gradually formed to neutralize the effective surface charge on the interface, $\mathrm{Cl^-}$ and $\mathrm{Na^+}$ away from the interface reach a dynamic equilibrium.
Physically, there is a big reservoir with inexhaustible $\mathrm{Cl^-}$ and $\mathrm{Na^+}$. 
When one looks at a ball $B_L=\{|\bm{x}|<L\}$ large enough, in order to ensure the conservation of density and momentum, the influx of ions should be equivalent to the outflux of ions through $\partial B_L$ in the sense of charges and heat. Although an ion would not change its direction at once when it crosses $\partial B_L$, there would be another ion from the reservoir which gets inside. Since we only care about the statistical behavior of $\mathrm{Cl^-}$ and $\mathrm{Na^+}$, we can simply bounce an ion back when it crosses $\partial B_L$ as if there was a virtual wall. This motivates us to consider an approximate problem in a truncated domain and impose a reflecting BC in the artificial wall. See section \ref{sec: approx} for mathematical justification of introducing a truncated domain for the PB equation in symmetric cases. 
 
\begin{figure}[H]
\centering
\includegraphics[width=0.3\textwidth]{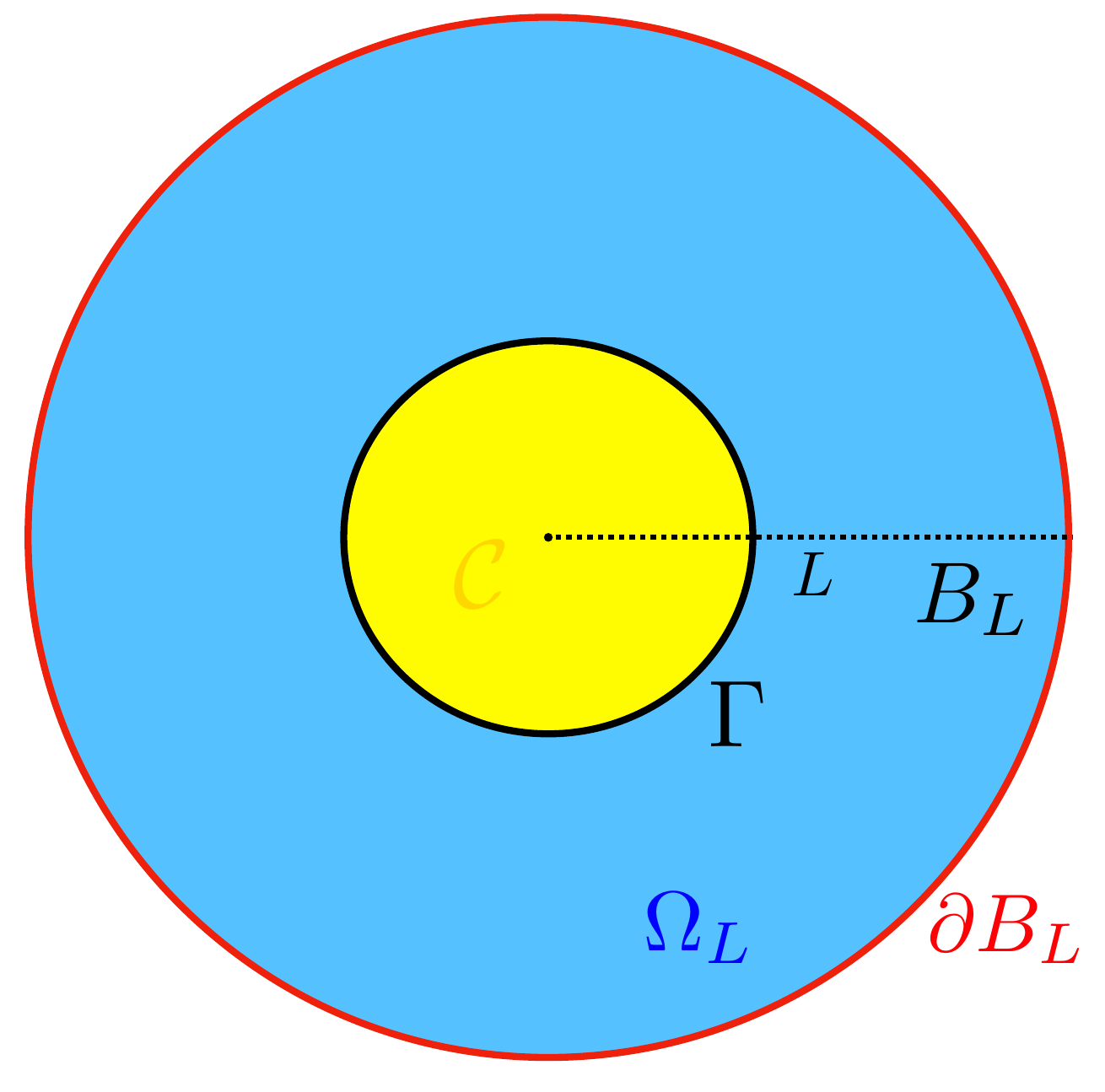}
\caption{Illustration of the truncated domain. The yellow and blue regions denote the object $\mathcal{C}$ and the truncated domain $\Omega_L$, while the black and red curves denote the inner boundary $\Gamma$ and the artificial boundary $\partial B_L$ respectively.}
\label{fig: illustr2}
\end{figure}

Let $B_L=\{|\bm{x}|<L\}$ be a sufficiently large ball containing $\mathcal{C}$.Then the truncated domain is $\Omega_L=B_L\backslash \bar{\mathcal{C}}$, see Figure \ref{fig: illustr2}. As mentioned before, one can only simulate the motion of charged particles in the truncated domain $\Omega_L$. Besides the reflecting process on the inner boundary $\Gamma$, we also use a reflecting BC for $\partial B_L$. Thus the approximate RSDE with artificial wall reads
\begin{equation}
\label{eq:RSDE}
\left\{
\begin{array}{ll}
	d X_j=-\nabla(z_j\hat{\Phi}_L)dt+\sqrt{2} dB+dR,\\[3mm]
    R(t)=\int_0^t\bm{n}(X_j)\mathrm{d}|R|_s,\quad |R|_t=\int_0^t\mathbbm{1}_{\partial \Omega_L}(X_j)\mathrm{d}|R|_s,\quad X_j(0)=X_{0,j}\sim \rho_{0,j}.
\end{array}
\right.
\end{equation} 
Here, the subscript $L$ is used to emphasize the truncated domain. The physical potential $\hat{\Phi}_L$ should be determined by Gauss's law, whose source consists of the free charges inside $\mathcal{C}$ and the positive and negative charges in $\Omega_L$. Moreover, $\partial \Omega_L=\Gamma\cup \partial B_L$, the reflecting process $R$ ensures the particle stay in $\Omega_L$. 

\subsection{The approximate PNP and PB equations with an artificial wall}

At the macroscopic level, the Nernst-Planck equations in $\Omega_L$ corresponding to \eqref{eq:RSDE} are
\begin{equation}
\label{eq: bounded PNP}
\left\{
\begin{array}{lll}
	\partial_t \rho_j= \nabla\cdot(z_j\rho_j\nabla\hat{\Phi}_L+\nabla\rho_j),\quad \bm{x} \in \Omega_L, \\[3mm]
	\langle z_j\rho_j\nabla\hat{\Phi}_L+\nabla\rho_j,\bm{n}\rangle=0,\quad \bm{x} \in \partial \Omega_L,\\[3mm]
	\rho_j|_{t=0}=\rho_{0,j},\quad \bm{x}\in \Omega_L.
\end{array}
	\right.\quad j=\pm.
\end{equation}
The no-flux BCs ensure that the total positive charge $Q_+=\int_{\Omega_L} \rho_+\mathrm{d}\bm{x}$ and total negative charge $Q_-=\int_{\Omega_L} \rho_-\mathrm{d}\bm{x}$ are conserved.

The total potential $\hat{\Phi}_L$ is generated by the free charges $\rho_f$ in $\mathcal{C}$ and the charges $\rho_{\pm}$ in $\Omega_L$. Due to the assumption of uniform dielectric constant, one has
\begin{equation}
\label{eq:physical potential}
	\hat{\Phi}_L=(\rho_f+\rho_+-\rho_-)*\Psi. 
\end{equation}
Here $\Psi$ is the Coulomb potential given by
\begin{equation}
\label{eq: Coulomb potential}
	\Psi(\bm{x})=\left\{
	\begin{array}{lll}
		-\frac{1}{2\nu} |\bm{x}|,&\quad d=1,\\[3mm]
-\frac{1}{2\pi\nu}\ln |\bm{x}|,&\quad d=2,\\[3mm]
\frac{1}{d(d-2)\alpha(d)\nu|\bm{x}|^{d-2}},&\quad d\geq 3,
	\end{array}
    \right.
\end{equation}
where $\alpha(d)=\frac{\pi^{d/2}}{\Gamma(d/2+1)}$ denotes the volume of the unit ball in $\mathbb{R}^d$. Physically, since the dielectric constants are the same in $\Omega$ and $\mathcal{C}$, there is no induced charge on the inner boundary $\Gamma$ as well as the artificial boundary $\partial B_L$, so $\Psi$ is the fundamental solution to $-\nu\Delta\Psi=\delta$ in the whole space. See Appendix \ref{sec: dielectric} for more discussions on variable dielectric constants.

The stationary solutions to \eqref{eq: bounded PNP} still have the form
$$\rho_{+}=C_{+}e^{-\hat{\Phi}_L},\quad \rho_{-}=C_{-}e^{\hat{\Phi}_L}.$$
Then the steady state Poisson equation reads
\begin{equation}
\label{eq: PB with C}
	-\nu \Delta \hat{\Phi}_L=C_{+}e^{-\hat{\Phi}_L}-C_{-}e^{\hat{\Phi}_L},\quad \bm{x}\in \Omega_L,
\end{equation}
$C_{+}, C_{-}$ are to be determined.

For the original PB equation \eqref{eq:PB}, $\rho_{+}^{\infty}=\rho_{-}^{\infty}=\rho_{\infty}$. Otherwise,
the net charge of the physical system $Q_f+\int_{\Omega}\rho_{+}-\rho_{-}\mathrm{d}\bm{x}$ is infinite since $\rho_+\rightarrow \rho_+^{\infty}, \rho_{-}\rightarrow \rho_{-}^{\infty}$ as $|\bm{x}| \rightarrow \infty$, where
$Q_f=\int_{\mathcal{C}} \rho_f \mathrm{d}\bm{x}$ is the total free charge in $\mathcal{C}$. 
However, in the bounded case, it is not necessary that $\rho_{+}=\rho_{-}$ for all $\bm{x} \in\partial B_L$. But we can impose
\begin{equation}
\label{eq: condition on rho}
	\rho_{+}(\bar{\bm{x}})\rho_{-}(\bar{\bm{x}})=\rho_{\infty}^2
\end{equation}
at some point $\bar{\bm{x}}\in \partial B_L$ by adjusting the total charge of positive ions. For example, in 1D, we can require $\rho_{+}(L)\rho_{-}(L)=\rho_{\infty}^2$. While in 3D, we can pick a point $\bar{\bm{x}}=(L,0,0)$ such that  $\rho_{+}(\bar{\bm{x}})\rho_{-}(\bar{\bm{x}})=\rho_{\infty}^2$. Thus, $C_{+}C_{-}=\rho_{\infty}^2$, and there exists a constant $c$ such that 
\begin{equation}
\label{eq: C}
	 C_{+}=\rho_{\infty}e^{-c},\quad C_{-}=\rho_{\infty}e^{c}.
\end{equation}
Let $\Phi_L=\hat{\Phi}_L+c$. We obtain 
\begin{equation}
\label{eq: bounded PB}
	\left\{
\begin{array}{llll}
-\nu \Delta \Phi_L=\rho_f ,\quad \bm{x}\in \mathcal{C},\\[1.5mm]
	-\nu \Delta \Phi_L=\rho_{\infty}\left(e^{-\Phi_L}-e^{\Phi_L}\right),\quad \bm{x}\in \Omega_L,\\[2mm]
		[\Phi_L]\big|_{\Gamma}=0, \quad [\nabla \Phi_L\cdot \bm{n}]\big|_{\Gamma}=0,\\[2.5mm]
	\dfrac{\partial \Phi_L}{\partial \bm{n}}\bigg|_{\partial B_L}=\bm{n}\cdot\nabla\left[ \left(\rho_f+\rho_{\infty}e^{-\Phi_L}-\rho_{\infty}e^{\Phi_L}\right)*\Psi\right],
\end{array}
\right.
\end{equation}
which is the approximate PB equation in the truncated domain.
We also remark that 
\begin{equation}
\label{eq: BI}
    \int_{\partial B_L}\frac{\partial \Phi_L}{\partial \bm{n}}\mathrm{d}S_{\bm{x}}= 0
\end{equation}
due to zero net charge at $\partial B_L$.

\section{Error estimate of the truncation for the PB equation}
\label{sec: approx}

In this section, we will show in the PDE level the validity of the truncation in a special case for the PB equation. In other words, the solution of the PB equation in the truncated domain $\Omega_L$ converges to the solution of the PB equation in $\Omega$ exponentially fast. We believe similar error control happens for the PNP equations but we leave this for future investigation.

For the sake of simplicity, we consider the case where $\rho_f$ is radially symmetric. Hence, 
$-\frac{\partial \Phi}{\partial \bm{n}}\big|_{\Gamma} =\sigma_f,$
where $\sigma_f$ denotes the equivalent effective surface charge for the inside free charge distribution $\rho_f$. This means that $\sigma_f$ is a given function totally determined by $\rho_f$. $\sigma_f$ gives exactly the same field as $\rho_f$ for $\bm{x}\in \Omega$ and it holds that
$ \nu\oint_{\Gamma}\sigma_f \mathrm{d} S_{\bm{x}}=Q_f$.
Then the PB equation in the external domain $\Omega$ can be rewritten as
\begin{equation}
    \label{eq:PB sym}
	\left\{
\begin{array}{ll}
	-\nu\Delta \Phi=\rho_{\infty}\left(e^{-\Phi}-e^{\Phi
	}\right),\quad \bm{x}\in \Omega,\\[2mm]
		-\dfrac{\partial \Phi}{\partial \bm{n}}\bigg|_{\Gamma} =\sigma_f,\quad \Phi(\bm{x})\rightarrow 0 \ \text{as} \ |\bm{x}| \rightarrow \infty.
	\end{array}
\right.
\end{equation}
For the truncated problem, zero electronic flux in $\Omega_L$ implies the electronic field is tangential. That is to say, \eqref{eq: BI}  will reduce to $\frac{\partial \Phi_L}{\partial \bm{n}}\big|_{\partial B_L}=0$ in this case. 
Then, the approximate problem in $\Omega_L=B_L\backslash \bar{\mathcal{C}}$ is 
\begin{equation}
\label{eq:Approx}
	\left\{
\begin{array}{ll}
	-\nu\Delta \Phi_L=\rho_{\infty}\left(e^{-\Phi_L}-e^{\Phi
	_L}\right),\quad \bm{x}\in \Omega_L,\\[2mm]
	-\dfrac{\partial \Phi_L}{\partial \bm{n}}\bigg|_{\Gamma} =\sigma_f,\quad \dfrac{\partial \Phi_L}{\partial \bm{n}}\bigg|_{\partial B_L} =0.
\end{array}
\right.
\end{equation}
The well-posedness of \eqref{eq:Approx} is guaranteed by classical elliptic theory \cite{Evans2010Partial}. By the way, in 1D, the constraint on $\rho_f$ can be removed due to the symmetry of the electronic field $-\nabla \Psi=\frac{1}{2\nu} \sgn(x)$. 

Next, we will show \eqref{eq:Approx} is a good approximation to \eqref{eq:PB sym} when the ball $B_L$ is large enough. The approximate result here focuses on 1D and 3D. However, we believe a similar convergence rate holds in other dimensions due to the exponential decay of $\nabla \Phi$. This can be estimated by analyzing the Green function of the linearized equation by the boundary integral method \cite{Hackbusch1995Integral,johnson2012numerical} like 3D. 

\begin{theorem}[\bf{Convergence}]
\label{thm: convergence}
$(d=1\& 3)$. Let $\Phi$ be the solution of \eqref{eq:PB sym} and $\Phi_L$ that of \eqref{eq:Approx}. Then there exist constants $C_1>0, C_2>0$ independent of $L$ such that
$$
\|\Phi_L-\Phi\|_{L^1(\Omega_L)}\leq  C_2 e^{-\frac{C_1}{\sqrt{\nu}} L}
$$
for $L$ large enough.
\end{theorem}

In order to show this, we need the exponentially decay property of the solution to \eqref{eq:PB sym}.

\begin{proposition}[\bf{Exponential decay of $\Phi, \nabla \Phi$}]
\label{prop: decay rate}
 $(d=1\& 3)$. There exists a positive constant $R$ and generic constants $C$ such that, for all $ |\bm{x}|>R$, the following estimates hold:
\begin{itemize}
	\item $d=1$.
\begin{align*}
	\left|\Phi(x)\right|\leq & ~ \frac{|\sigma_f|}{\sqrt{2\rho_{\infty}/\nu}}e^{-\sqrt{2\rho_{\infty}/\nu}\rm{dist}(x,\Gamma)},\\[2mm]
	\left|\Phi^{\prime}(x)\right|\leq & ~ |\sigma_f| e^{-\sqrt{2\rho_{\infty}/\nu}\rm{dist}(x,\Gamma)}.
\end{align*}

	\item $d=3$.
	\begin{align*}
	\left|\Phi(\bm{x})\right|\leq & ~ \frac{C}{|\bm{x}|}e^{-\frac{C}{\sqrt{\nu}}|\bm{x}|},\\[2mm]
	\left|\nabla\Phi(\bm{x})\right|\leq & ~ \frac{C}{|\bm{x}|^2}e^{-\frac{C}{\sqrt{\nu}}|\bm{x}|}.
\end{align*}
\end{itemize}
\end{proposition}
This proposition is proved by comparing it to a linear equation, whose solution is given by the boundary integral representation formula. So we need the following preparations: Lemma \ref{lem: representation formula} and Lemma \ref{lem: comparison principle}. The proof of Proposition \ref{prop: decay rate} in 3D is left to Appendix \ref{sec: proofs}. For 1D, it is quite similar and simpler, so it is omitted.

\begin{lemma}[\cite{johnson2012numerical} \bf{Representation formula}]
\label{lem: representation formula}
If $u$ is smooth in $\mathbb{R}^d\backslash \Gamma$ and satisfies 
$$
\left\{
\begin{array}{ll}
-\Delta u+cu=f \quad \text{in}\ \Omega,\\[2mm]
-\Delta u+cu=0 \quad \text{in}\ \mathcal{C},\\[2mm]
\lim\limits_{|\bm{x}| \rightarrow \infty}u=0,\quad \lim\limits_{|\bm{x}| \rightarrow \infty}|\nabla u|=0.
\end{array}\right.
$$
Then
\begin{equation}
\label{eq: representation formula}
	\begin{aligned}
&\int_{\Gamma}\left\{\left[\frac{\partial u}{\partial \bm{n}}\right]G(\bm{x}-\bm{y})-[u]\frac{\partial G(\bm{x}-\bm{y})}{\partial \bm{n}_{\bm{y}}}\right\}\mathrm{d}S_{\bm{y}}+\int_{\Omega}G(\bm{x}-\bm{y})f(\bm{y})\mathrm{d}\bm{y}\\[2mm]
=&
\left\{
\begin{array}{ll}
	\vspace{2mm} u(\bm{x}),&\quad \bm{x}\notin \Gamma,\\
	\dfrac{u^i(\bm{x})+u^e(\bm{x})}{2},&\quad \bm{x}\in \Gamma.
\end{array}
\right.
\end{aligned}
\end{equation}
Here, $G$ is the Green function which is the solution to 
$(cI-\Delta) G =\delta$. For $ \bm{x}\in \Gamma$, $u^i(\bm{x})$ and $u^e(\bm{x})$ represent the limit from $\mathcal{C}$ and $\Omega$, respectively. $[u]=u^i-u^e$.
\end{lemma}
The proof is similar to that in \cite{johnson2012numerical} and we omit the details.

\begin{lemma}[\cite{Evans2010Partial} \bf{Comparison principle on unbounded domain}]
\label{lem: comparison principle}
Let $u\in C^2(\Omega)\cap C(\bar{\Omega})$.
Consider an elliptic operator $\mathcal{L}$ having the form
$$\mathcal{L}u=-\Delta u + c(u) u,$$
where $c>0$ is continuous. Then
\begin{equation*}
	\left\{
\begin{array}{lll}
    \vspace{2mm}\mathcal{L}u\leq 0,\quad \text{in} \ \Omega,\\
	\vspace{2mm}-\dfrac{\partial u}{\partial \bm{n}}\leq 0,\quad \text{on} \ \Gamma,\\
	\lim\limits_{|\bm{x}|\rightarrow \infty} u(\bm{x})=0.
\end{array}
\quad \Rightarrow \quad u\leq 0.
\right.
\end{equation*}
\end{lemma}
We also sketch a short proof of Lemma \ref{lem: comparison principle} in Appendix \ref{sec: proofs}.

Now, we are ready to prove the convergence result: Theorem \ref{thm: convergence}. We only give a 3D proof, the 1D case can be shown without the bridge of sup solution.
\begin{proof}[Proof of Theorem \ref{thm: convergence} in 3D]
Construct a sup solution $\Phi^+$ which satisfies
\begin{equation}
\label{eq:sup}
	\left\{
\begin{array}{ll}
	-\nu\Delta \Phi^+=\rho_{\infty}\left(e^{-\Phi^+}-e^{\Phi^+}\right),\quad \bm{x}\in \Omega_L,\\[3mm]
	-\dfrac{\partial \Phi^+}{\partial \bm{n}}\Big|_{\Gamma} =\sigma_f,\quad \dfrac{\partial \Phi^+}{\partial \bm{n}}\Big|_{\partial B_L} =\Sigma_L,
\end{array}
\right.
\end{equation}
where $\Sigma_L=\left\|\dfrac{\partial \Phi}{\partial \bm{n} }\right\|_{L^{\infty}(\partial B_L)}< +\infty$.

Let $u_1=\Phi^+ -\Phi_L$. According to \eqref{eq:sup} and \eqref{eq:Approx}, $u_1$ satisfies 
\begin{equation}
\label{u1}
	 \left\{
\begin{array}{ll}
	-\nu \Delta u_1+cu_1=0,\quad \bm{x} \in \Omega_L,\\[3mm]
	-\dfrac{\partial u_1}{\partial \bm{n}}\Big|_{\Gamma} =0,\quad \dfrac{\partial u_1}{\partial \bm{n}}\Big|_{\partial B_L} =\Sigma_L,
\end{array}
\right.
\end{equation}
with
$c=\dfrac{2\rho_{\infty}\left(\sinh \Phi^+ -\sinh\Phi_L\right)}{\Phi^+ -\Phi_L}=2\rho_{\infty}\cosh\Phi_{\xi_{u_1}}\geq 2\rho_{\infty}$ bounded below on $\Omega_L$. Then by the maximum principle, we obtain 
\begin{equation}
	u_1\geq 0. 
\end{equation}

Integrating equation \eqref{u1} on $\Omega_L$ yields
$$\nu \int_{\partial B_L}\Sigma_L \mathrm{d}S_{\bm{x}}=\nu \int_{\partial \Omega_L}\frac{\partial u_1}{\partial \bm{n}} \mathrm{d}S_{\bm{x}}=\int_{\Omega_L} c u_1\mathrm{d}\bm{x}.$$
It follows from the non-negativity of $u_1$ that
$$\|u_1\|_{L^1(\Omega_L)}=\left|\int_{\Omega_L} u_1 \mathrm{d}\bm{x}\right|\leq \frac{2\pi L^2\nu}{\rho_{\infty}}\Sigma_L.$$
That is, 
$$\|\Phi^+ -\Phi_L\|_{L^1(\Omega_L)}\leq \frac{2\pi L^2\nu }{\rho_{\infty}}\|\nabla \Phi\cdot \bm{n}\|_{L^{\infty}\left(\partial B_L\right)}
\leq \frac{2\pi L^2\nu }{\rho_{\infty}}\|\nabla \Phi\|_{L^{\infty}\left(\partial B_L\right)} .$$

Similarly, let $u_2=\Phi^+ -\Phi$ for $\bm{x} \in \bar{\Omega}_L$.
Then, $$\|\Phi^+ -\Phi\|_{L^1(\Omega_L)}\leq \frac{4\pi L^2\nu }{\rho_{\infty}}\|\nabla \Phi\|_{L^{\infty}\left(\partial B_L\right)} .$$
Hence, one concludes from Proposition \ref{prop: decay rate} that there exist positive constants $R$ and $C_1, C_2$ such that for $L>R$,
$$\|\Phi_L-\Phi\|_{L^1(\Omega_L)}\leq  C_2 e^{-\frac{C_1}{\sqrt{\nu}}L}.$$
\end{proof}
Theorem \ref{thm: convergence} demonstrates the exponentially decay rate of $L^1$ error of $\Phi$ and $\Phi_L$ in truncated domain $\Omega_L$. This verifies we only need to solve the PB system in $\Omega_L$.

\section{The random batch particle methods}
\label{sec: RBM}

In \eqref{eq:RSDE}, the law of particles is the positive or negative charge density distribution up to a multiplicative constant.
Hence, \eqref{eq:RSDE} can be regarded as the mean-field limit of the interacting particle system. In this section, we investigate the interacting particle system and state three essential problems in solving it. We then propose some numerical methods to obtain our random batch particle methods.

\subsection{Interacting particle systems for the PNP and PB equations}
\label{subsec: system}

In \eqref{eq:RSDE}, the potential is the one generated by the mean-field distribution. We use $N$ particles that interact with each other through Coulomb interaction to approximate the distribution, hoping that the empirical measures of positive and negative particles $\rho_{\pm}^N$ multiplied by total positive and negative charge $Q_{\pm}$ are the approximations to $\rho_{\pm}$. That is, in the $N\rightarrow \infty$ regime,
$$
Q_{+}\rho_{+}^N\rightharpoonup \rho_{+},\quad Q_{-}\rho_{-}^N\rightharpoonup \rho_{-},
$$
where$$
\rho_{+}^N=\frac{1}{N_+}\sum_{i\in I_+}\delta(\cdot-X^i),\quad \rho_{-}^N=\frac{1}{N_{-}}\sum_{i\in I_{-}}\delta(\cdot-X^i).
$$
Here, the superscript $i$ denotes the $i$-th particle, $I_+=\{i, z^i=1\}, I_{-}=\{i, z^i=-1\}$ with $z^i$ being the sign of the $i$-th particle. Let $q$ be the absolute value of charge per particle. Then, the numbers of positive and negative charges are $N_+=Q_+/q$ and $N_-=Q_-/q$ respectively. The total particle number is $N=N_{+}+N_{-}$. Denote $|Q|=Q_++Q_-$ as the total absolute charge in $\Omega_L$, we have $|Q|=Nq$. Note that the 'particle' here can be either 'numerical particle' or 'physical particle'.  

Let $$F=-\nabla \Psi,\quad E_f=\rho_f*F,$$ so that
$F=\frac{\bm{x}}{d\alpha(d)\nu|\bm{x}|^d}$ is the Coulomb repulsive force and $E_f$ is the electronic field generated by $\rho_f$.
Then, the above interpretation implies that we can approximate the self-consistent RSDE \eqref{eq:RSDE} through the interacting particle system
\begin{equation}
\label{eq: RSDE_NP}
\left\{
\begin{aligned}
	dX^i=& ~z^iE_f(X^i)dt+\sum_{k:k\neq i}z^iz^kqF(X^i-X^k)dt+\sqrt{2} dB^i+dR^i,\ \  i=1,\cdots,N,\\
	R^i(t)=&\int_0^t\bm{n}(X^i)\mathrm{d}|R^i|_s,\quad |R^i|_t=\int_0^t\mathbbm{1}_{\partial \Omega_L}(X^i)\mathrm{d}|R^i|_s,\quad X^i(0)=X^i_0.
\end{aligned}
\right.
\end{equation}
Here $\{X^i\}_{i=1}^{N}$ are the trajectories of $N$ particles and the particles from all species are numbered together. $\{B^i\}_{i=1}^{N}$ are $N$ independent $d$-dimensional Brownian motions and $\{R^i\}_{i=1}^{N}$ are reflecting processes associated with $\{X^i\}_{i=1}^N$. The initial data $\{X^i_0\}_{i\in I_{\pm}}$ are independent, identically distributed (i.i.d.) random variables with probability density function $\rho_{0,\pm}$. 

Next, we point out three essential issues in solving the $N$-particle system \eqref{eq: RSDE_NP}.

The first one is unphysical attraction. Intuitively, as $N\to\infty$, \eqref{eq: RSDE_NP} will approximate \eqref{eq:RSDE}, and this is the so-called mean-field limit. In fact, if the force field $F$ is regular, such mean-field limit can be justified rigorously.
For the interacting particle system \eqref{eq: RSDE_NP} with finite $N$ particles, however, there is some probability that two opposite particles attract each other. When a positive and a negative charge meet, they cancel each other and there is an energy jump (the interaction energy is frozen and set to zero). This is a problem related to the $N$-particle system \eqref{eq: RSDE_NP}, which does not arise in the $N\to\infty$ limit. In fact, in the $N\rightarrow \infty$ limit, each particle carries infinitely small charge and this cancellation will not affect the continuum interaction energy. 
	
A possible way to deal with this issue is to introduce hard sphere potential as in the physical model with finite number of ions. The hard sphere potential due to finite size effect will avoid the unphysical attraction and cancellation. Likewise, we can add hard sphere potentials in the $N$-particle system \eqref{eq: RSDE_NP}. The Lennard-Jones potential \cite{LJ1924}
	\begin{equation}
	   	\phi(\bm{x})=4\epsilon\left[\left(\frac{\sigma}{|\bm{x}|}\right)^{12}-\left(\frac{\sigma}{|\bm{x}|}\right)^6\right],
	\label{eq: LJ} 
	\end{equation}
where $\epsilon$ is the depth of the potential well and $\sigma$ is the finite distance at which the inter-particle potential is zero, is usually used. 
As a consequence, the interacting particle system becomes 
	\begin{equation}
	    dX^i=z^i E_f(X^i)dt+\sum_{k:k\neq i}z^iz^kqF(X^i-X^k)dt-\sum_{k:k\neq i}\nabla \phi(X^i-X^k)dt+\sqrt{2} dB^i+dR^i.
	    \label{eq: LJ system}
	\end{equation}
In the mean-field limit $N\to\infty$, the parameters $\epsilon$ and $\sigma$ will vanish and the mean-field limit of \eqref{eq: LJ system} is then \eqref{eq:RSDE}. Hence, \eqref{eq: LJ system} can be used for the $N$-particle computation. Meanwhile, as the hard sphere potential does not play its role in the mean-field limit, we do not need the hard sphere potential if we aim to capture the mean-field behaviors only, as long as we can find a way to tackle the attraction issue. We will see this as a byproduct of RBM in section \ref{subsubsec: time split RBM}. 

The second issue for \eqref{eq: RSDE_NP} or \eqref{eq: LJ system} is the singularity of the forces which brings in numerical stiffness. We will discuss how to resolve it using time splitting strategy in section \ref{subsubsec: time split RBM} or kernel splitting strategy in section \ref{subsubsec: kernel split RBM}.

What matters most is the computational cost. Direct simulation of \eqref{eq: RSDE_NP} or \eqref{eq: LJ system} is expensive due to the interaction term. Solving it requires $\mathcal{O}(N^2)$ operation per time step.
As discussed in the introduction, there are several strategies to resolve this, and some typical methods include the FMM or the PIC. 
The FMM and PIC could have better accuracy with $\cO(N)$ or $\cO(N\log N)$ cost, but the prefactor in the linear scaling could be large. 
We choose to apply the random batch method (RBM) in \cite{JinLiLiu20} to \eqref{eq: RSDE_NP} or \eqref{eq: LJ system}, which is also $\mathcal{O}(N)$ but has smaller prefactor in the linear scaling. The reasons include the simplicity for implementation and better scalability in parallel computing (as the random-batch based RBE method demonstrates in \cite{liang2021superscalability}). Moreover, compared to PIC, such particle methods allow the incorporation of more physical effects. 
Note that the random batch method is based on the Monte Carlo ideas, so the accuracy is not very high if the step size is not very small and batch size is not very big. However, when the accuracy requirement is not very high, it can still potentially have much less CPU cost (\cite{JinLiLiu20,JinLiSun20,RBE}). Also, RBM will introduce additional noise due to the randomness but this can be controlled by the interaction with the heat bath in the Langevin equations.  The comparison between the RBM strategy and the traditional methods has been thoroughly discussed in \cite{JinLiLiu20,JinLiSun20,RBE,liang2021superscalability}.

Below, we briefly introduce the RBM strategy and then discuss how to implement the RBM strategy for our systems to turn into practical random batch particle methods in section \ref{subsec: numer method}.

\subsection{Random batch particle methods for the PNP and PB equations}
\label{subsec: numer method}

\subsubsection{A brief introduction to the RBM}
\label{subsubsec: RBM}

Pick a time step $\tau>0$ and define the time grid $t_m:=m\tau$. On each time sub-interval $[t_{m-1}, t_m)$, RBM \cite{JinLiLiu20,JinLiSun20} randomly divides the $N$ particles into $n$ small batches with batch size $p$ ($p\ll N$) and interact them within each batch. For the next time step, one reshuffles and forms a new set of batches, and repeats the process. Note that the set of random batches can be obtained in $\mathcal{O}(N)$ cost using the random permutation so that the computational cost per time step is remarkably reduced from $\mathcal{O}(N^2)$ to $\mathcal{O}(pN)$ per time step. The idea used to fasten the evaluation of interacting force shares a lot of similarities with the Direct Simulation Monte Carlo method \cite{Bird70,Nanbu80} based on binary collisions for Boltzmann equation and its adaptation for mean-field equations of flocking dynamics using stochastic binary interactions \cite{AlPa13}. RBM has been investigated theoretically in \cite{JinLi21,JinLiLiuweights20}, and already has a variety of applications in, for example, efficient sampling \cite{SVGDRBM,RBMC}, 
MD \cite{RBE,liang2021superscalability}, flocking models \cite{HaJinKimKo} and quantum systems \cite{GJP,QMC20}.

The original version of RBM in \cite{JinLiLiu20} was for indistinguishable particles, but it was then extended to interacting particles with disparate species and weights \cite{JinLiLiuweights20}. In our case, the particles have different charges so we will apply the version in \cite{JinLiLiuweights20}, which will be explained briefly here. Consider the first order interacting particle system
\begin{equation}
\label{eq: full particle system}
dX^i = b(X^i)dt +\frac{1}{N-1}\sum_{k:k\neq i}m_kF_{ik}(X^i,X^k)dt+\sigma dB^i,\ \ i=1,\cdots,N.
\end{equation}
 For each time interval $[t_{m-1},t_m)$, RBM solves the following SDE instead
\begin{equation}
\label{eq: RBM-1}
dX^i = b(X^i)dt +\frac{1}{p-1}\sum_{k\in \mathcal{C}_{\theta(i)}^{(m)}:k\neq i}m_kF_{ik}(X^i,X^k)dt+\sigma dB^i.
\end{equation}
Here, $\mathcal{C}^{(m)}:=\{\mathcal{C}_{\ell}^{(m)}: \ell=1,\cdots, n\}$ denotes the random batches on $[t_{m-1}, t_m)$, and the batches will be renewed at next time grid point. $\theta(i)$ indicates the index $\ell$ such that $i\in \mathcal{C}_{\ell}^{(m)}$. The cost is clearly $\cO(pN)$ per time step for the new system since interactions only take place inside small batches.

Let us briefly explain why RBM works here. Define the fluctuation of the random force on particle $i$ by 
\begin{equation}
\chi_i= \frac{1}{p-1}\sum_{k\in \mathcal{C}_{\theta(i)}:k\neq i}m_k F_{ik}(\bm{x}^i,\bm{x}^k)-\frac{1}{N-1}\sum_{k:k\neq i}m_k F_{ik}(\bm{x}^i,\bm{x}^k).
\end{equation}
It is proved in {\cite[Lemma 3.2]{JinLiLiuweights20}} that 
\begin{equation}\label{eq:fluc}
    \mathbb{E}\chi_i =0, \quad \mathrm{Var}(\chi_i)=\mathbb{E}|\chi_i|^2=\left(\frac{1}{p-1}-\frac{1}{N-1}\right)\Lambda_i,
\end{equation}
where 
\[
\Lambda_i=\frac{1}{N-2}
        \sum_{j: j\neq i}\Big| m_j F_{ij}(\bm{x}^i, \bm{x}^j)-\frac{1}{N-1}
        \sum_{\ell: \ell\neq i}m_{\ell} F_{i\ell}(\bm{x}^i, \bm{x}^{\ell})  \Big|^2
\]
and the expectation is taken over random divisions for a given configuration $(\bm{x}^1,\cdots, \bm{x}^N)$. $\Lambda_i$ is independent of batch size $p$.
This claims the random force is unbiased and the variance is smaller for larger $p$. For a single time step, the cross-batch interactions are completely neglected, which leads to an $\cO(1)$ approximation error (i.e. $\chi_i=\cO(1)$). However, since we do random reshuffling at each time grid, the random errors will roughly cancel out over time as the random force is unbiased by \eqref{eq:fluc}. As the dynamics go on, the time averaging effect owing to the law of large numbers (in time) could ensure the convergence of RBM. Since the number of time intervals is like $\tau^{-1}$ so the strong error would be like $\sqrt{\mathrm{Variance}*\tau}\sim \sqrt{\tau/p}$ by a typical Monte Carlo bound.

Since the configurations on the time intervals are not independent, the error bound cannot be obtained from the law of large numbers directly. Nevertheless, it has been shown in \cite{JinLiLiuweights20} that the strong error of RBM for regular interactions can be given by the following result. 
\begin{proposition}[{\cite[Theorem 3.1]{JinLiLiuweights20}}]
\label{prop: strong convergence}
Let $X^i$ and $\tilde{X^i}$ be solutions to \eqref{eq: full particle system} and \eqref{eq: RBM-1} respectively. Suppose the weights $m_i$ are bounded, the external force $b$ is one-sided Lipschitz, $b, \nabla b$ have polynomial growth, and the interacting force $F_{ik}$ have uniformly bounded second order derivatives. Then there exists $C$ independent of $N, p$ such that
\begin{equation}
\label{eq: strong error RBM} 
		\mathop{\sup}_{t\leq T}J(t):=\mathop{\sup}_{t\leq T}\sqrt{\frac{1}{2N}\sum_{i=1}^Nm_i\mathbb{E}\left|X^i(t)-\tilde{X}^i(t)\right|^2 }\leq  C\sqrt{\frac{\Lambda}{p-1}\tau},
	\end{equation}
	where $\Lambda=\max_i ||\Lambda_i||_{\infty}$.
\end{proposition}
The error bound clearly agrees with the Monte Carlo interpretation above. Hence, the strong convergence order is $1/2$ and larger batch size $p$ gives more accurate approximation.
A more important implication of the error estimate is that the error bound is uniform in $N$ so that the method is ``asymptotic-preserving" under the mean-field limit (see \cite{JinLiLiu20} for more discussions). Hence, we can choose batch size $p=\cO(1)$ independent of $N$ and the method really scales like $\cO(N)$.

If one cares about the distribution generated by RBM and the statistics of the particle system (like density, pressure, etc), the weak error makes more sense. It has been shown in \cite{JinLiLiuweights20} that the distribution generated by RBM is indeed close to the one generated by the full particle system and the error in the weak sense is first order.
\begin{proposition}[{\cite[Theorem 4.1]{JinLiLiuweights20}}]
\label{prop: weak convergence}
Suppose the weights $m_i$ are bounded, the functions $b, F_{ik}$ are $C^4$ and have uniformly bounded derivatives up to order $4$. Then, for any test function $\varphi\in C_b^{\infty}(\mathbb{R}^d)$, the weak error is controlled as
	\begin{equation}
\label{eq: weak error RBM} 
\mathop{\sup}_{m:t_m\leq T}E_m:=\mathop{\sup}_{m:t_m\leq T}\left|\frac{1}{N}\sum_{i=1}^N \omega_i \mathbb{E}\varphi(\tilde{X}^i(t_m))-\frac{1}{N}\sum_{i=1}^N \omega_i \mathbb{E}\varphi(X^i(t_m))\right|\leq C\tau,
\end{equation}
where $ \omega_i=\frac{N m_i}{\sum_{k=1}^N m_k}$, $C=C(\varphi,T)$ is independent of $N, \tau$.
\end{proposition}
Estimate \eqref{eq: weak error RBM} shows that the RBM system \eqref{eq: RBM-1} converges weakly with first order to the full system \eqref{eq: full particle system}, in terms of empirical measures. Standard mean-field theory \cite{JabinWang17,JabinWang18} says the empirical measure of $N$-particle system converges to the solution to the corresponding Fokker-Planck equation as $N\rightarrow \infty$ and the invariant measure as $t\rightarrow \infty$. Hence, we expect that we can use RBM to approximate the distributions of positive and negative charges in the PNP and PB equations, though the rigorous proof is still open since the Coulomb interaction is singular.

\subsubsection{An effective particle method for the PNP and PB equations}
\label{subsubsec: time split RBM}

We now consider the RBM approximation to \eqref{eq: RSDE_NP} to reduce the cost.
The corresponding RBM system reads: for $t\in[t_{m-1},t_m)$ and $i=1,\cdots, N$, 
\begin{eqnarray}
\label{eq: RBM-PNP}
\left\{
\begin{aligned}
	dX^i=&~z^iE_f(X^i)dt+\frac{N-1}{p-1}\sum\limits_{k\in \mathcal{C}_{\theta(i)}^{(m)}:k\neq i}z^iz^kqF(X^i-X^k)dt\\
	&+\sqrt{2} dB^i+dR^i,\\
	R^i(t)=&\int_0^t\bm{n}(X^i)\mathrm{d}|R^i|_s,\quad |R^i|_t=\int_0^t\mathbbm{1}_{\partial \Omega_L}(X^i)\mathrm{d}|R^i|_s,\quad X^i(0)=X^i_0.
\end{aligned}
\right.
\end{eqnarray}

Unlike the full particle system \eqref{eq: RSDE_NP}, in random batch particle system \eqref{eq: RBM-PNP}, if two opposite particles encounter (they are in the same batch and the distance between them is zero) at a certain time step, it is of high probability that they get lost (they are in different batches) in the next time step. Owing to the random mini batch approximation, two particles being stuck all the time is an impossible event, which does not have to be taken care of. This means the attraction issue automatically disappears when using RBM. Hence, system \eqref{eq: RBM-PNP} is an effective interacting particle system for the mean-field behavior described by the PNP and PB equations. Discretization of this will yield an effective particle method for the PNP and PB equations.

The stiffness in \eqref{eq: RBM-PNP} due to the singularity can be resolved well in the case $p=2$. In fact, using the time splitting method, we may split system \eqref{eq: RBM-PNP} into 
\begin{align}
&dX^i=~(N-1)z^iz^kqF(X^i-X^k)dt, \quad i,k \in \mathcal{C}_{\ell}^{(m)}, \label{eq: split 1}\\
&dX^i=~z^iE_f(X^i)dt+\sqrt{2} dB^i+dR^i\label{eq: split 2}.
\end{align}

Therefore, for each time step, we can solve \eqref{eq: split 1} analytically and then apply stochastic schemes to \eqref{eq: split 2}. 
This then gives an effective random batch particle method for mean-field behavior described by the PNP and PB equations, as detailed in Algorithm \ref{algo: RBM PNP/PB}.

\begin{algorithm}[H]
\caption{An effective random batch particle method for the PNP and PB equations}
\begin{algorithmic}[1]
\For{ $m$ in $1:[T/\tau]$} 
\State Divide $\{1,2,\cdots,N\}$ into $n=N/2$ batches randomly.
\For{each batch $\mathcal{C}_{\ell}^{(m)}$}
\State Update $X^i,X^k (i,k\in\mathcal{C}_{\ell}^{(m)})$ for $t\in [t_{m-1},t_m)$ by the following:
\State Compute 
$v=\frac{X^i-X^k}{\abs{X^i-X^k}},\quad \beta=\frac{2(N-1)z^iz^kq}{\alpha(d)\nu },\quad \eta=\abs{X^i-X^k}^d+\beta \tau $.
\State If $\eta \geq 0$,
 $\tilde{X}^i=\frac{X^i+X^k+v\eta^{1/d}}{2},\quad
    \tilde{X}^k=\frac{X^i+X^k-v\eta^{1/d}}{2}.$
\State Otherwise, $\tilde{X}^i=\tilde{X}^k=\frac{X^i+X^k}{2}.$
\State Solve \eqref{eq: split 2} for $X^i,X^k$ with initial data $\tilde{X}^i,\tilde{X}^k$.
\EndFor
\EndFor
\end{algorithmic}
\label{algo: RBM PNP/PB}
\end{algorithm}

\subsubsection{The kernel splitting strategy for general batch sizes and molecular dynamics}
\label{subsubsec: kernel split RBM}

The time splitting method above clearly does not work for \eqref{eq: RBM-PNP} when $p\ge 3$. Instead, we consider the splitting strategy introduced in \cite{split,MTMC} and decompose $F$ into two parts: 
$$
F=F_1+F_2,
$$
where $F_1$ is regular with long-range, $F_2$ is singular with short-range. 
In principle, the splitting is kind of arbitrary as long as the stiff part is included in the short-range part. A feasible splitting for the Coulomb force when $d\geq 2$ could be 
\begin{equation}
    F_1=\left\{
\begin{array}{cc}
    f_c,& \quad |\bm{x}|<r_c,\\
    \frac{\bm{x}}{d\alpha(d)\nu |\bm{x}|^d} ,&\quad |\bm{x}|>r_c,
\end{array}
\right.
\quad
F_2=\left\{
\begin{array}{cc}
     \frac{\bm{x}}{d\alpha(d)\nu |\bm{x}|^d}-f_c,&  \quad |\bm{x}|<r_c,\\
     0 ,& \quad |\bm{x}|>r_c,
\end{array}
\right.
\label{eq: decomposition}
\end{equation}
where $r_c$ is the cutoff radius, $f_c=\frac{r_c}{d\alpha(d)\nu |r_c|^d}$.
Applying RBM to $F_1$ only yields the split RBM approximation to \eqref{eq: RSDE_NP}. With suitable $r_c$, the summation in $F_2$ can be done in $\mathcal{O}(1)$ operation for each particle $X^i$. So, the overall cost is still $\mathcal{O}(N)$ per time step. 
Clearly, this kernel splitting approach for general batch sizes $p\ge 3$ can be used and the resulting RBM system reads
\begin{equation}
\begin{aligned}
	dX^i=& ~\frac{N-1}{p-1}\sum\limits_{k\in \mathcal{C}_{\theta(i)}^{(m)}:k\neq i}z^iz^kqF_1(X^i-X^k)dt
	+\sum_{k:k\neq i}z^iz^kqF_2(X^i-X^k)dt\\
	&+z^iE_f(X^i)dt+\sqrt{2} dB^i+dR^i.
\end{aligned}
\label{eq: RBM_split}
\end{equation}

Since the short-range attraction is computed fully in \eqref{eq: RBM_split}, the attraction issue may arise. Hence, one would like to include the hard sphere potentials. Moreover, if one considers the physical models, the Lennard-Jones potential \eqref{eq: LJ} would be essential so one would like to simulate \eqref{eq: LJ system}.
The kernel splitting strategy mentioned above can be used for random batch particle methods corresponding to \eqref{eq: LJ system}, either for MD simulations or for numerical method of the PNP and PB equations. The resulting random batch particle method is shown in Algorithm \ref{algo: split RBM}.
Here, the decomposition \eqref{eq: decomposition} has been applied for the Coulomb potential while Lennard-Jones potential has not been split due to the short-range nature. Of course, if one desires, the Lennard-Jones potential may also be suitably decomposed in applications. When used as a numerical method, our experience shows that such kind of method is comparable to Algorithm \ref{algo: RBM PNP/PB} for solving PNP and PB equations, see section \ref{subsec: sym 3D exmp}.

\begin{algorithm}[H]
\caption{Split RBM for the PNP and PB equations}
\begin{algorithmic}[1]
\For{ $m$ in $1:[T/\tau]$} 
\State Divide $\{1,2,\cdots,N=pn\}$ into $n$ batches randomly.
\For{each batch $\mathcal{C}_{\ell}^{(m)}$}
\State Update $X^i$'s $(i\in\mathcal{C}_{\ell}^{(m)})$ for $t\in [t_{m-1},t_m)$ by the following
\begin{equation}
\begin{aligned}
	dX^i=& ~\frac{N-1}{p-1}\sum\limits_{k\in \mathcal{C}_{\ell}^{(m)}:k\neq i}z^iz^kqF_1(X^i-X^k)dt
	+\sum_{k:k\neq i}z^iz^kqF_2(X^i-X^k)dt\\
	&-\sum_{k:k\neq i}\nabla \phi(X^i-X^k)dt+z^iE_f(X^i)dt+\sqrt{2} dB^i+dR^i.
\end{aligned}
\label{eq: RBM_LJ}
\end{equation}
\EndFor
\EndFor
\end{algorithmic}
\label{algo: split RBM}
\end{algorithm}

\subsection{An iterative RBM-PB method with fixed $\rho_{\infty}$ } 

In many applications, we only care about the equilibrium and what we often know is the far field concentration $\rho_{\infty}$ in solution, rather than the total positive or negative charge $Q_{\pm}$. However, in the above methods, $Q_{\pm}$ is assumed to be known so that we know how many positive or negative particles we need. So, we design an iterative algorithm to determine $Q_{+}$ with fixed $\rho_{\infty}$.

\begin{figure}[H]
\centering
\includegraphics[width=0.6\textwidth]{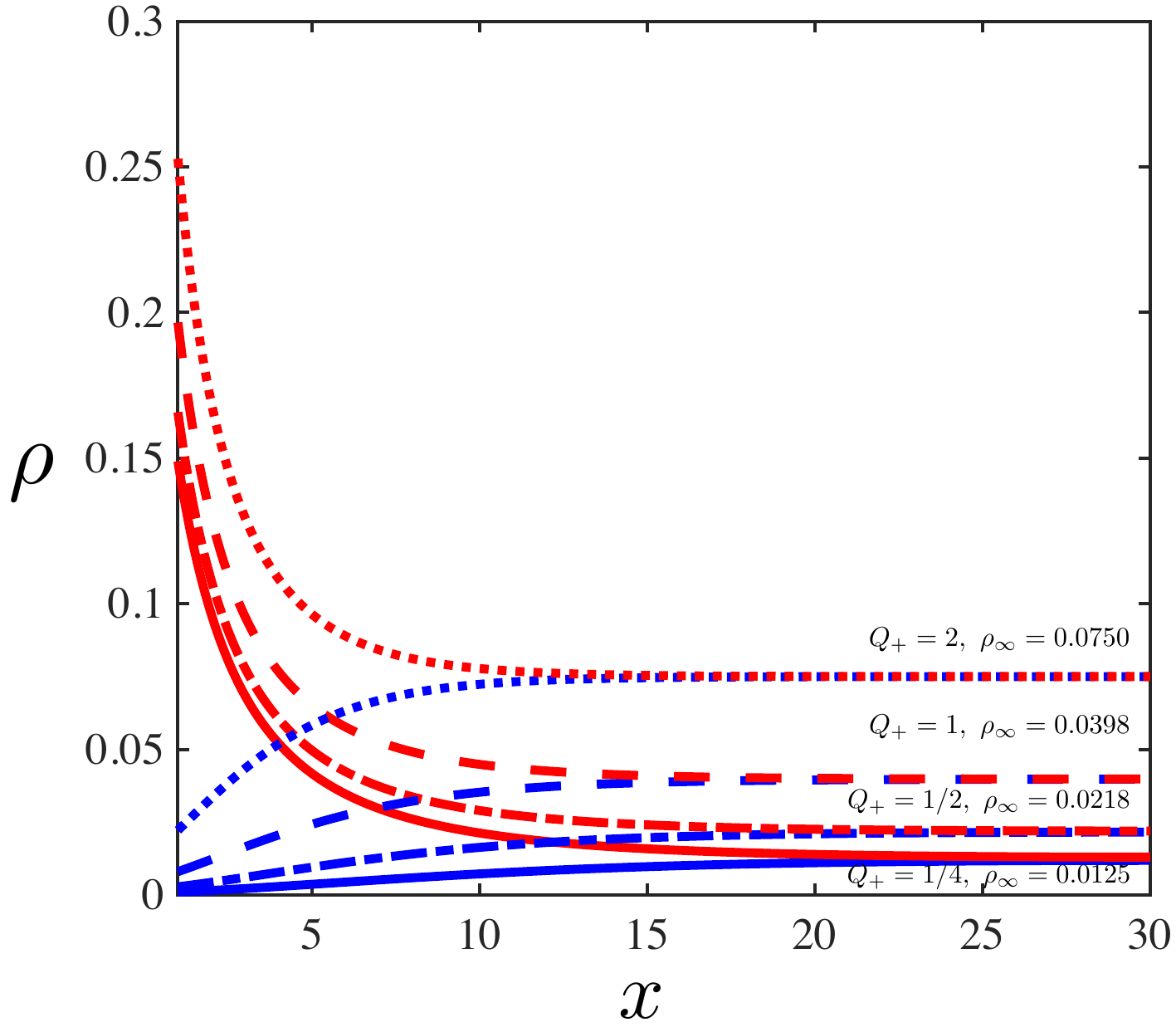}
\caption{The equilibrium distributions (red for $\rho_-$, blue for $\rho_+$) in 1D with fixed $Q_f=1$. The truncated domain is $[1, 10]$. This shows $\rho_{\infty}$ is larger when $Q_+$ is larger.}
\label{fig:fact}
\end{figure}

One can see from Figure \ref{fig:fact} that: with fixed $Q_f$, the larger $Q_+$ is, the larger $\rho_{\infty}$ will be. Though plotted in 1D, it holds for higher dimensions. Since each particle shares charge $q$, larger $Q_+$ represents more positive particles. So we can run the particle simulation by iteration without prior knowledge of the total positive charge $Q_+$ and adjust the particles adaptively. 

In each iteration, first run RBM simulation for the PB equation with $Q_+$ at present. Then compute the densities $\rho_{\pm}(\bar{\bm{x}})$ at some fixed point $\bar{\bm{x}}\in \partial B_L$. If $\rho_+(\bar{\bm{x}})\rho_-(\bar{\bm{x}})$ is lower than $\rho_{\infty}^2$, randomly add some positive and negative particles into the ensemble. Otherwise, randomly kill some positive and negative particles. Since the free charge distribution $\rho_f$ in $\mathcal{C}$ is unchanged, the total net charge $Q_{+} -Q_{-}$ in $\Omega_L$ should be fixed as $-Q_f$, which leads to a synchronous change of positive and negative particles. The iteration terminates when condition \eqref{eq: condition on rho} is satisfied within tolerance. After iteration, one can get the true $Q_+$. The process is illustrated in Algorithm \ref{algo: iterative RBM}. 
 
\begin{algorithm}[H]
\caption{Iterative RBM-PB method with fixed $\rho_{\infty}$}
{\bf Input} Initial distributions $\rho_{0,\pm}$, truncated length $L$, free charge distribution $\rho_f$, charge per particle $q$, far field concentration $\rho_{\infty}$, a point $\bar{\bm{x}}\in \partial B_L$, tolerance $\delta$.
\begin{algorithmic}[1]
	\State Initial guess $Q_+$, $Err=1$.
	\State Generate $N_+=\frac{Q_+}{q}$ samples from $\rho_{0,+}$ and $N_-=\frac{Q_{-}}{q}=\frac{Q_{+}+Q_f}{q}$ samples from $\rho_{0,-}$.
	\While{$Err>\delta$}
	\State Run RBM simulation for PB with $Q_+$.
    \State Calculate the densities $\rho_{\pm}(\bar{\bm{x}})$. 
    \State $I=\rho_{+}(\bar{\bm{x}})\rho_{-}(\bar{\bm{x}})-\rho_{\infty}^2$, error $Err=\sqrt{|I|}$, the approximate charge change $\triangle Q:=\frac{\alpha(d) L^d}{2}Err$. Then $\triangle N:=\left\lfloor\frac{\triangle Q}{q}\right\rfloor, \triangle Q'=q\triangle N$.
    \If{$I<0$}
     \State Generate $\triangle N$ samples from $\rho_{0,\pm}$ and add them to the positive particles and the negative particles respectively. Set $Q_+\leftarrow Q_+ + \triangle Q'$.
     \Else
     \State Remove $\triangle N$ samples from the positive particles and negative particles respectively. Set $Q_+\leftarrow Q_+ - \triangle Q'$.
    \EndIf 
	\EndWhile
\end{algorithmic}
{\bf Output} $Q_+,\rho_+,\rho_-,\Phi_L$.
\label{algo: iterative RBM}
\end{algorithm}

In Algorithm \ref{algo: iterative RBM}, the bulk densities $\rho_{\pm}(\bar{\bm{x}})$ are
$$\rho_{\pm}(\bar{\bm{x}})=\frac{2N_{\pm}(D_h)}{\alpha(d)h^dN_{\pm}},$$ 
where $D_h$ is a small half ball $\{\bm{x}\in \Omega_L: |\bm{x}-\bar{\bm{x}}| \leq h\}$ around $\bar{\bm{x}}$, $N_{\pm}(D_h)$ is the number of positive or negative particles in $D_h$. 
Due to random fluctuation, the calculation of bulk densities  is not stable, thus the iteration may not converge. A direct way to avoid random fluctuation is to increase particle number, but this is expensive. Another widely used approach is time average, which is the average of densities of previous iterations after reaching equilibrium.

An alternative way to adjust particles adaptively is to take the grand canonical ensemble \cite{grandensemble80} into account. This will be explored in the future.

\section{Numerical examples}
\label{sec: numerics}

In this section, we give a 1D example and a symmetric 3D example to show that the random batch particle methods can be successfully used as the numerical methods for solving the PNP and PB equations. Then, a non-symmetric 3D example and a colloidal example are given to illustrate the adaptivity to complicated geometry and complex phenomenon respectively.

\subsection{1D case}
\label{subsec: 1D exmp}
Consider a 1D example with $\mathcal{C}=(-R, R)$. Assume the free charge distribution inside $\mathcal{C}$ is $\rho_f=\sum_{i=1}^Kq_i\delta(x-x_i), x_i\sim U(-R, R)$, then the total free charge $Q_f=\sum_{i=1}^Kq_i$. Outside $\mathcal{C}$, the solution is symmetric 1:1 salt. As mentioned before, we truncate the external domain and only focus on the positive half $\Omega_L=(R, L)$ due to symmetry. We use $ -\partial_x\Phi_L(R)=\sigma_f$, $\partial_x\Phi_L(L)=0$ 
as the BCs for the approximate PNP and PB equations. 
By the superposition principle, the effective surface charge for the right half $\sigma_f=E_f(R)+E_r(R)=\sum_{i=1}^K\frac{q_i}{2}\frac{1}{2\nu}-\frac{Q_f}{2}(-\frac{1}{2\nu})=\frac{Q_f}{2\nu}$, where $E_f=\sum_{i=1}^K\frac{q_i}{2}\frac{1}{2\nu}\sgn(x-x_i)$ is the external field for the right half domain, $E_r$ is the field generated by particles in solution with total net charge $-\frac{Q_f}{2}$. The initial data for the PNP equations are chosen as uniform distributions in $(R, L)$ with $c_{\pm}=Q_{\pm}/(L-R)$ for simplicity. We comment that the initial distributions do not affect the equilibrium but will affect the dynamics. Besides, $\Phi_L$ is unique up to a constant for the PNP equations with Neumann BCs, we further impose $\mathbb{E}(\Phi_L)=0$ in computation.

Note that in the simulations, the particles can cross each other. This can be regarded as a special case in 3D, where the distribution is homogeneous in the $y, z$ directions so that only the $x$ direction is left. Particles are thus charge sheets in 3D so they can cross each other freely. 

\begin{figure}[H]
\centering
\includegraphics[width=13.5cm,height=5cm]{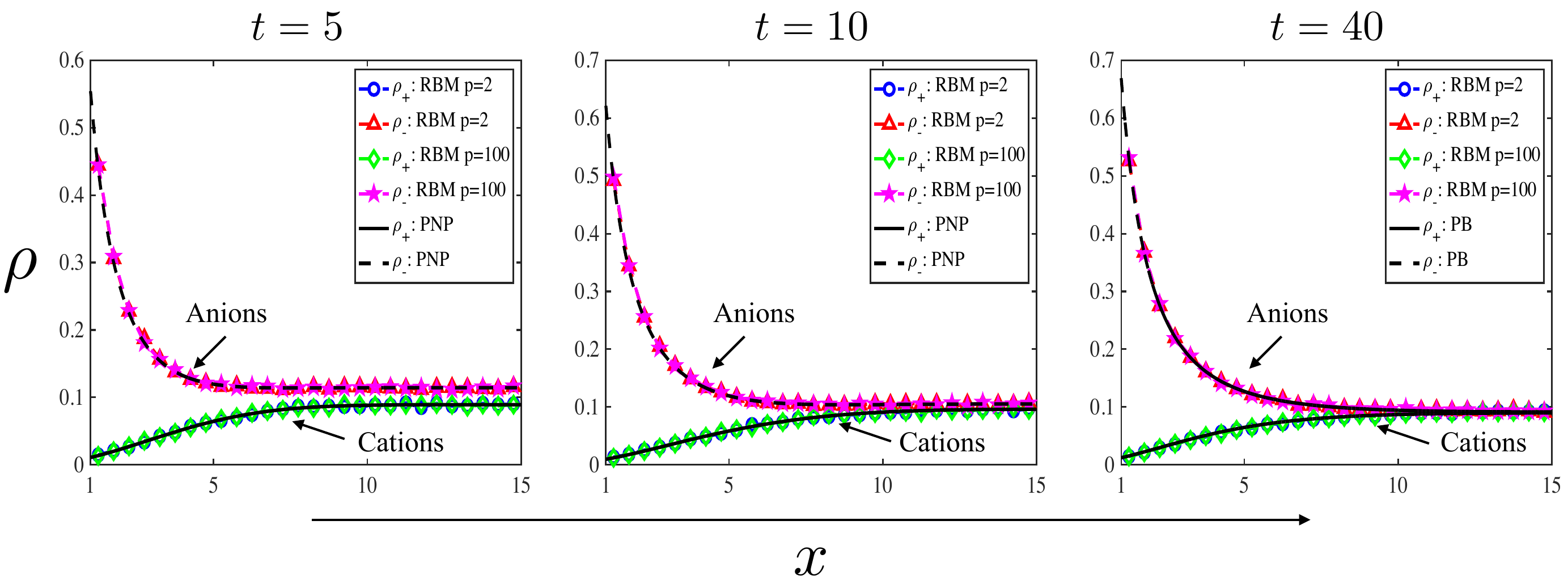}
\caption{Comparison of RBM results with PNP at $t=5, 10$ and with PB at $t=40$.}
\label{fig:1D dyn}
\end{figure}

Let $\nu=1, \rho_f=\sum_{i=1}^{100}q_i\delta (x-x_i)$ with $q_i\sim U(-3,3)$, $R=1, L=15$, $Q_{+}=1$, $q=1e-5$. We show the performance of RBM with batch size $p=2$ and $p=100$ and step size $\tau=0.01$. Since the 1D Coulomb force $F=\frac{1}{2\nu}\sgn (x)$ is regular, no splitting strategy is needed. For the dynamic problem, the reference solution is given by the conservative FD scheme of the PNP equations.
While for the equilibrium, the reference solution is given by Newton's iteration using FD to the PB equation. The comparison of distributions at different times is in Figure \ref{fig:1D dyn}. We can observe that the RBM simulation results match well with the reference solutions to both PNP and PB equations.

\begin{figure}[H]
\centering
\includegraphics[width=14cm,height=5cm]{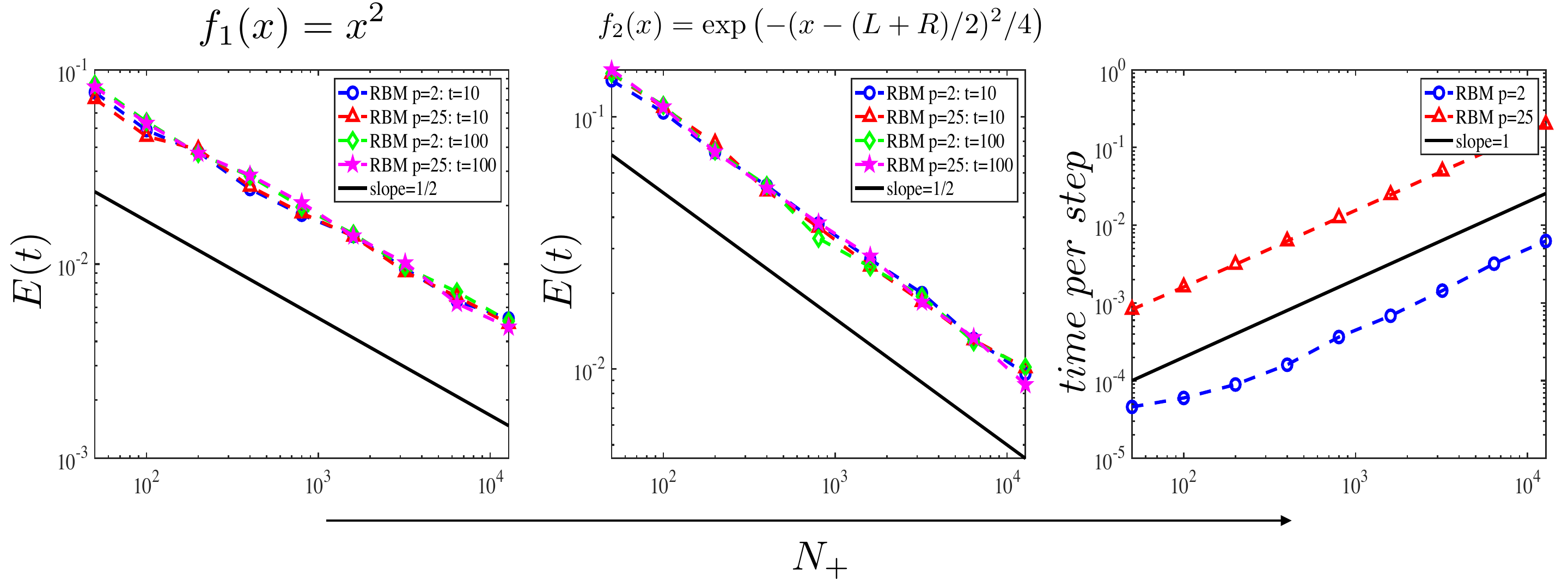}	
\caption{The weak error of RBM corresponding to PNP at $t=10$ (left) and PB at $t=100$ (middle) and the cost of RBM(right) versus positive particle number $N_{+}$.}
\label{fig: 1D conv}
\end{figure}

Next, we test the weak convergence of RBM. For a given test function $f(x)$, we compute the estimated expectation for positive and negative charges
$\bar{f}_{\pm}(t)=\frac{1}{N_{\pm}}\sum_{i\in I_{\pm}} f(X^i(t))$. 
Define the Mean Square Error (MSE) over $M$ independent repetitive experiments
$$
\mathrm{MSE}_{\pm}(t):=\frac{1}{M}\sum_{m=1}^M\left(\frac{\bar{f}^{(m)}_{\pm}(t)-\mathbb{E}_{x\sim \rho^*_{\pm}(t)}f(x)}{\mathbb{E}_{x\sim \rho^*_{\pm}(t)}f(x)}\right)^2.
$$
Here, $m$ means the $m$-th experiment, $\rho_{\pm}^*(t)$ are the reference densities of the PNP equation, $\mathbb{E}_{x\sim \rho^*_{\pm}(t)}f(x)=\int_{\Omega_L} f(x)\rho_{\pm}^*(t)dx$. The weak error is measured via the square root of average of the MSE of positive and negative charges
$$
E(t)=\sqrt{\frac{MSE_{+}(t)+MSE_{-}(t)}{2}}.
$$
We use two test functions $f_1(x)=x^2, f_2(x)=\exp\left(-(x-(L+R)/2)^2/4\right)$. Let $\nu=1, R=1, \rho_f=2\delta(x), Q_{+}=1, L=15, \tau=0.01, N_+=50, 100,\cdots, 12800, M=100$, the convergence results are shown in Figure \ref{fig: 1D conv}. Note that at equilibrium, we also have the weak error with $\rho^*_{\pm}$ being the solutions to the PB equation.
It is obvious that the RBM methods with both $p=2$ and $p=25$ are halfth-order in particle number $N$, which coincides with the Monte Carlo convergence rate. Besides, we plot the CPU time per time step in the last column of Figure \ref{fig: 1D conv}, which indicates $\mathcal{O}(N)$ computational cost clearly. 

\begin{figure}[H]
\centering
 \includegraphics[width=0.95\textwidth]{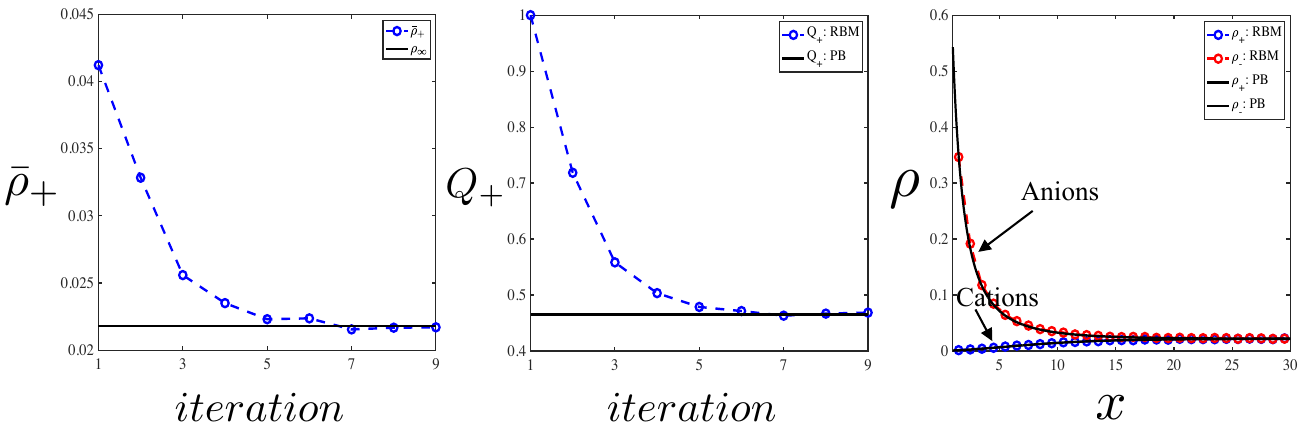}	
\caption{Left: trajectory of $\bar{\rho}_{+}$; middle: trajectory of $Q_{+}$; right: equilibrium distributions. }
\label{fig: 1D iteration}
\end{figure}

In former tests, we check the validity and convergence rate of RBM for PNP and PB equations with the total positive charge $Q_+$ given. Next, focus on the stationary problem and consider the problem where $\rho_{\infty}$ is given instead of $Q_+$ and test the iterative RBM-PB method proposed in Algorithm \ref{algo: iterative RBM}. Given $\rho_{\infty}=0.0218$, fix $\nu=1, R=1, L=30, \rho_f=\half\delta(x-\half)+\frac{3}{2}\delta(x+\half), q=1e-4, \delta=1e-5, \tau=0.1, p=2$, one test result is shown in Figure \ref{fig: 1D iteration}. After $9$ iterations, the RBM simulation reaches the equilibrium state in comparison with the reference result. This verifies that the efficient RBM method can also be applied to the case where $Q_+$ is unknown.

\begin{figure}[H]
\centering
\includegraphics[width=0.8\textwidth]{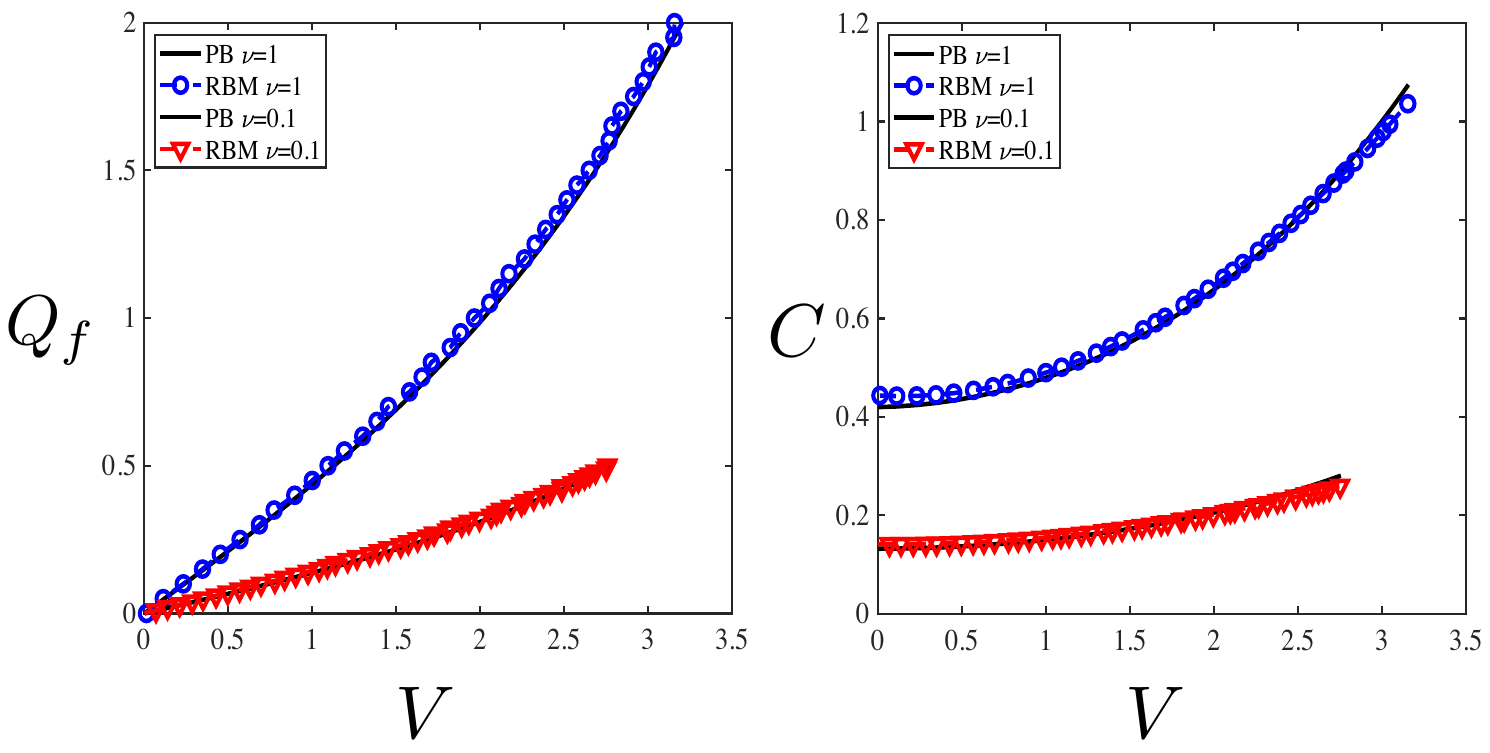}
	\caption{Differential capacitance for 1D example with fixed $\rho_{\infty}=0.0218$.}
	\label{fig: capacitance}
\end{figure}

Moreover, this physical model (a charged object in some electrolyte solution) can be regarded as a capacitor. The differential capacitance is defined by
$$C=\frac{d Q_f}{dV},$$
where $V=\Phi(R)-\Phi(L)$ is the voltage.
Figure \ref{fig: capacitance} shows it is feasible to measure the differential capacitance by iterative RBM-PB method under different scalings. The concentration of the electrolyte solution, i.e. $\rho_{\infty}$, is set to be $0.0218$. 

\subsection{3D spherical symmetric case}
\label{subsec: sym 3D exmp}

Assume the $\mathcal{C}$ is a ball centered at the origin with a free charge $Q_f$ at the center, i.e. $\mathcal{C} = \{\bm{x}\in\mathbb{R}^3: |\bm{x}|<R\}, \rho_f=Q_f\delta(\bm{x})$. As stated before, we simulate the particles inside a sufficiently large domain $\Omega_L=\{\bm{x}\in\mathbb{R}^3: R\leq |\bm{x}|\leq L\}$. Hence, the distributions of positive and negative charges in $\Omega$ (outside $\mathcal{C}$) are spherical symmetric. Denote $r=|x|$, consider the approximate PNP equations in the radial axis
\begin{equation}
    \label{eq: 3d pnp}
    \left\{
    \begin{array}{lllllll}
    \partial_t \rho_{+}=\partial_r\left(e^{-\Phi_L}\partial_r(e^{\Phi_L}\rho_{+})\right)+\frac{2}{r}e^{-\Phi_L}\partial_r(e^{\Phi_L}\rho_{+}),\quad r\in (R, L),\\[2mm]
    \partial_t \rho_{-}=\partial_r\left(e^{\Phi_L}\partial_r(e^{-\Phi_L}\rho_{-})\right)+\frac{2}{r}e^{\Phi_L}\partial_r(e^{-\Phi_L}\rho_{-}),\quad r\in (R, L),\\[2mm]
    e^{-\Phi_L}\partial_r(e^{\Phi_L}\rho_{+})=0,\quad r=R,L,\\[2mm]
    e^{\Phi_L}\partial_r(e^{-\Phi_L}\rho_{-})=0,\quad r=R,L,\\[2mm]
    \rho_{+}(r,0)=\frac{Q_{+}}{4\pi(L-R)r^2},\quad \rho_{-}(r,0)=\frac{Q_{-}}{4\pi(L-R)r^2},\quad r\in (R, L),\\[2mm]
    -\nu\left(\partial_{rr}\Phi_L+\frac{2}{r}\partial_r\Phi_L\right)=\rho_{+}-\rho_{-},\quad r\in (R, L),\\[2mm]
    -\partial_r\Phi_L(R)=\sigma_f,\quad \partial_r\Phi_L(L)=0
    \end{array}\right.
\end{equation}
and the approximate PB equation in the radial axis
\begin{equation}
    \label{eq: 3d pb}
     \left\{
    \begin{array}{ll}
    -\nu\left(\partial_{rr}\Phi_L+\frac{2}{r}\partial_r\Phi_L\right)=\rho_{\infty}(e^{-\Phi_L}-e^{\Phi_L}),\quad  r\in (R, L),\\[2mm]
    -\partial_r\Phi_L(R)=\sigma_f,\quad \partial_r\Phi_L(L)=0.
    \end{array}\right.
\end{equation}
Here, $\sigma_f=\frac{1}{4\pi R^2\nu}Q_f$ due to symmetry.

\begin{figure}[H]
\centering
\includegraphics[width=13.5cm,height=5cm]{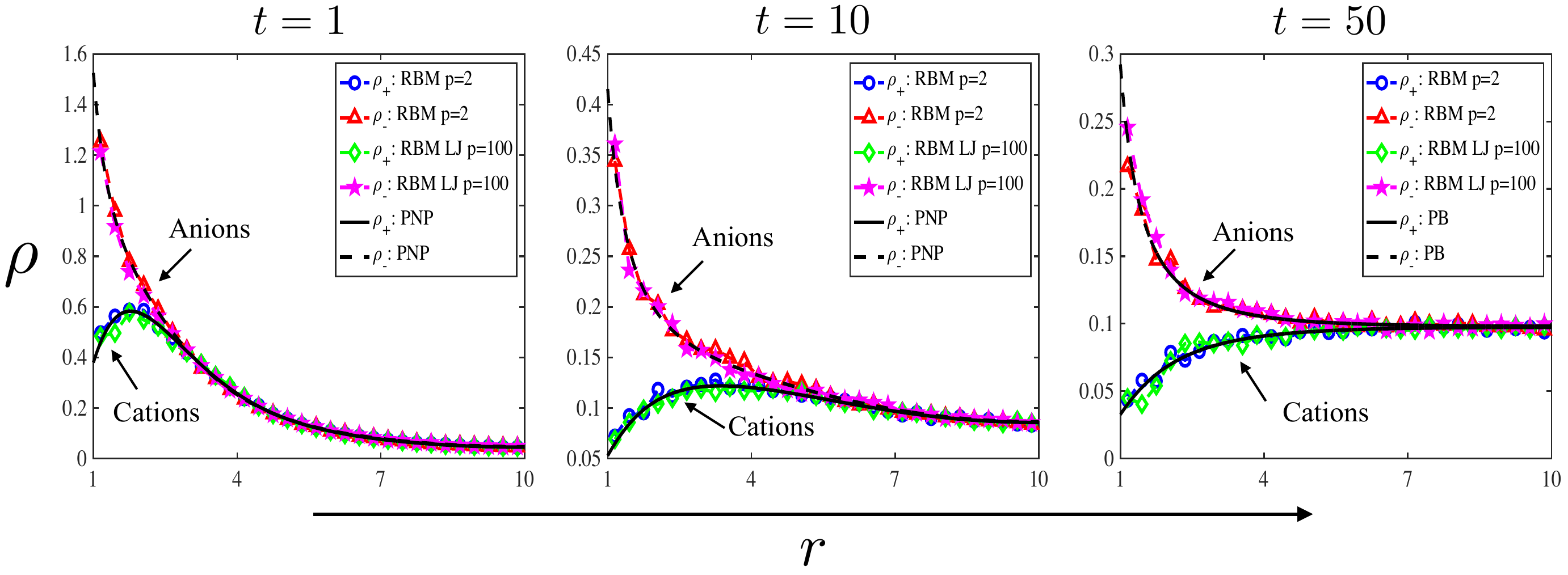}
\caption{Comparison of RBM results with PNP at $t=1, 10$ and with PB at $t=50$.}
\label{fig:3D dyn}
\end{figure}

We adopt the splitting strategies introduced in section \ref{subsec: numer method} to deal with the singular Coulomb interaction. Here, we do RBM simulations with batch size $p=2$ using time splitting in Algorithm \ref{algo: RBM PNP/PB} and RBM simulations including Lennard-Jones potential with batch size $p=100$ using kernel splitting in Algorithm \ref{algo: split RBM}. The reference solutions are also given by the conservative FD scheme to the PNP equations \eqref{eq: 3d pnp} and Newton's iteration to the PB equation \eqref{eq: 3d pb}. Take $\nu=1, R=1, L=10, Q_f=20\nu, Q_{+}=400$, in RBM simulations $q=4e-3, \tau=0.01$, the cut-off radius is $r_c=0.05$, the parameters in the Lennard-Jones are chosen as $\epsilon=1e-6q, \sigma=0.01$. The results are shown in Figure \ref{fig:3D dyn}. Again, we can conclude that RBM with proper splitting strategies can be used as numerical methods for both time-dependent PNP and stationary PB equations in 3D. Also, we can see that the simulation results with or without Lennard-Jones potential are comparable. This verifies our discussions that the use of random batch idea can avoid unphysical attraction in section \ref{sec: RBM}.

\begin{figure}[H]
\centering
\includegraphics[width=14cm,height=5cm]{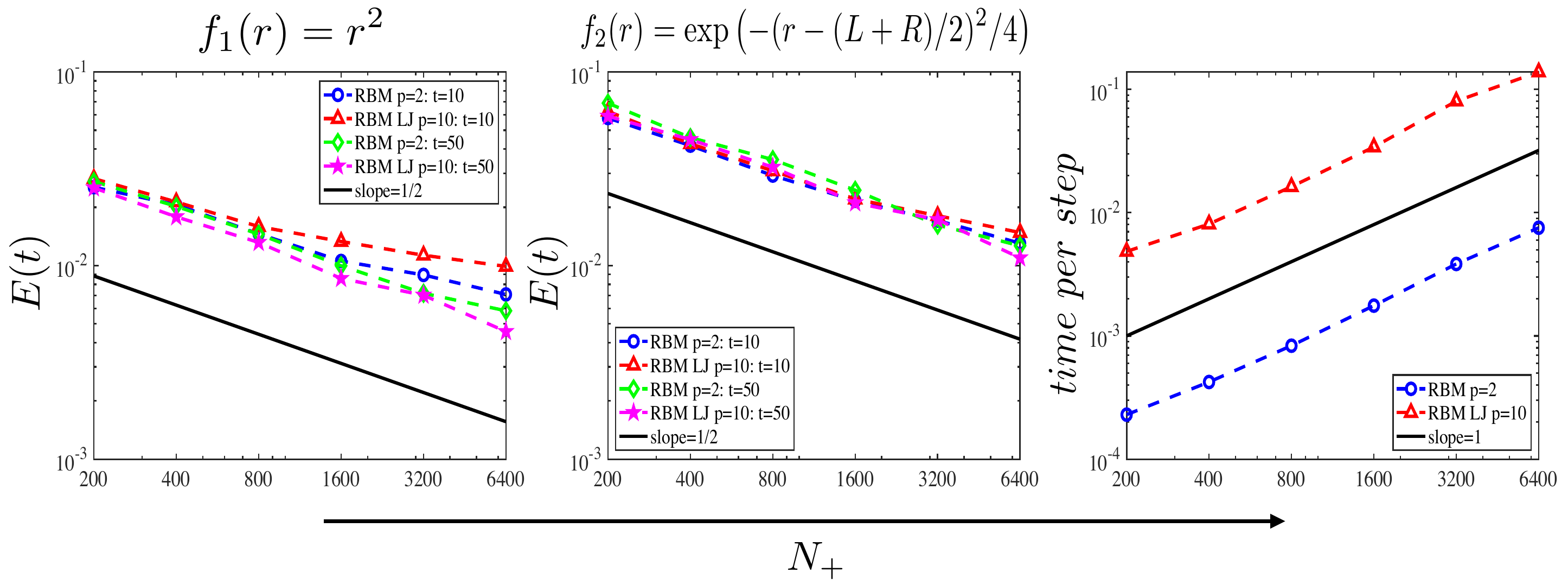}	
\caption{The weak error of RBM corresponding to PNP at $t=10$ (left), to PB at $t=50$ (middle) and the cost of RBM (right) versus positive particle number $N_{+}$. }
\label{fig: 3D conv}
\end{figure}

As in 1D case, we test the weak convergence of RBM. Let $\nu=1, R=1, Q_f=20\nu, Q_{+}=400, L=10, \tau=0.01, N_+=200,400,\cdots, 6400, r_c=0.05, \sigma=0.01, \epsilon=1e-6q, M=100$, the convergence results are shown in Figure \ref{fig: 3D conv}. Similarly, the RBM methods with both $p=2$ and $p=10$ including Lennard-Jones potential are halfth-order in particle number $N$, the cost is $\mathcal{O}(N)$. In terms of computational time, RBM with $p=2$ is more efficient, so we prefer time splitting RBM method (Algorithm \ref{algo: RBM PNP/PB}) in real computation.

\begin{figure}[H]
\centering
\includegraphics[width=0.8\textwidth]{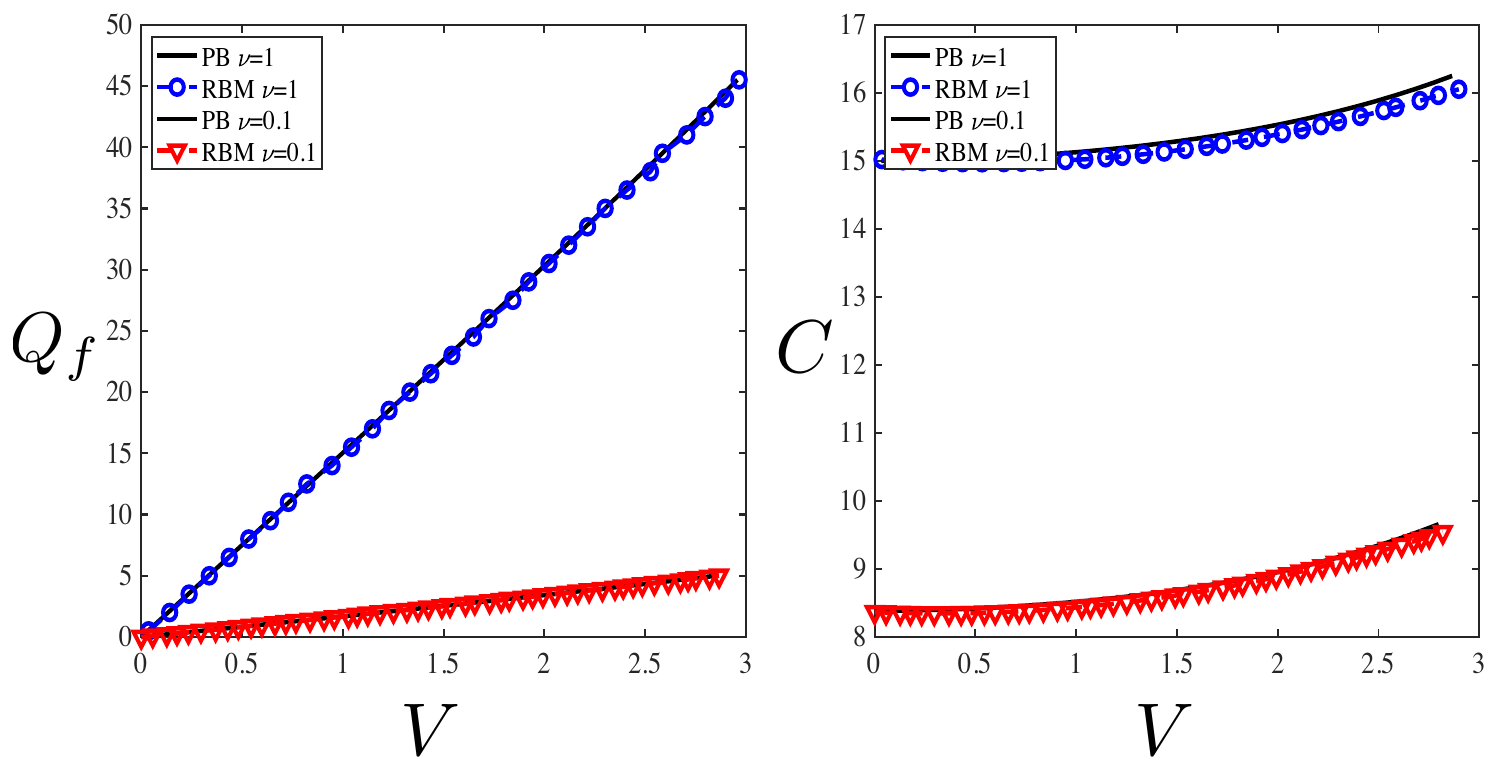}
	\caption{Differential capacitance for 3D spherical symmetric example with fixed $\rho_{\infty}=0.005$.}
	\label{fig: 3d capacitance}
\end{figure}

Similarly, given the ionic concentration $\rho_{\infty}$ of the electrolyte solution, one can compute the differential capacitance using the iterative RBM-PB method (Algorithm \ref{algo: iterative RBM}) in different scalings. We take $\rho_{\infty}=0.005, R=1, L=10, q=1e-3, \tau=0.01, \epsilon=1e-5$ and use the time splitting RBM to sample from the equilibrium. 
To weaken the random fluctuation, we collect the samples from $100$ time steps after the system reaches the equilibrium state. Therefore, we can use relatively small number of particles. 
The results for $\nu=1$ and $0.1$ are shown in Figure \ref{fig: 3d capacitance}. It is clear that the differential capacitance can be well approximated by RBM simulation for 3D spherical symmetric case. 

\subsection{3D non-symmetric case}
\label{subsec: non-sym 3D exmp}

\begin{figure}[H]
\centering
\includegraphics[width=5cm,height=4cm]{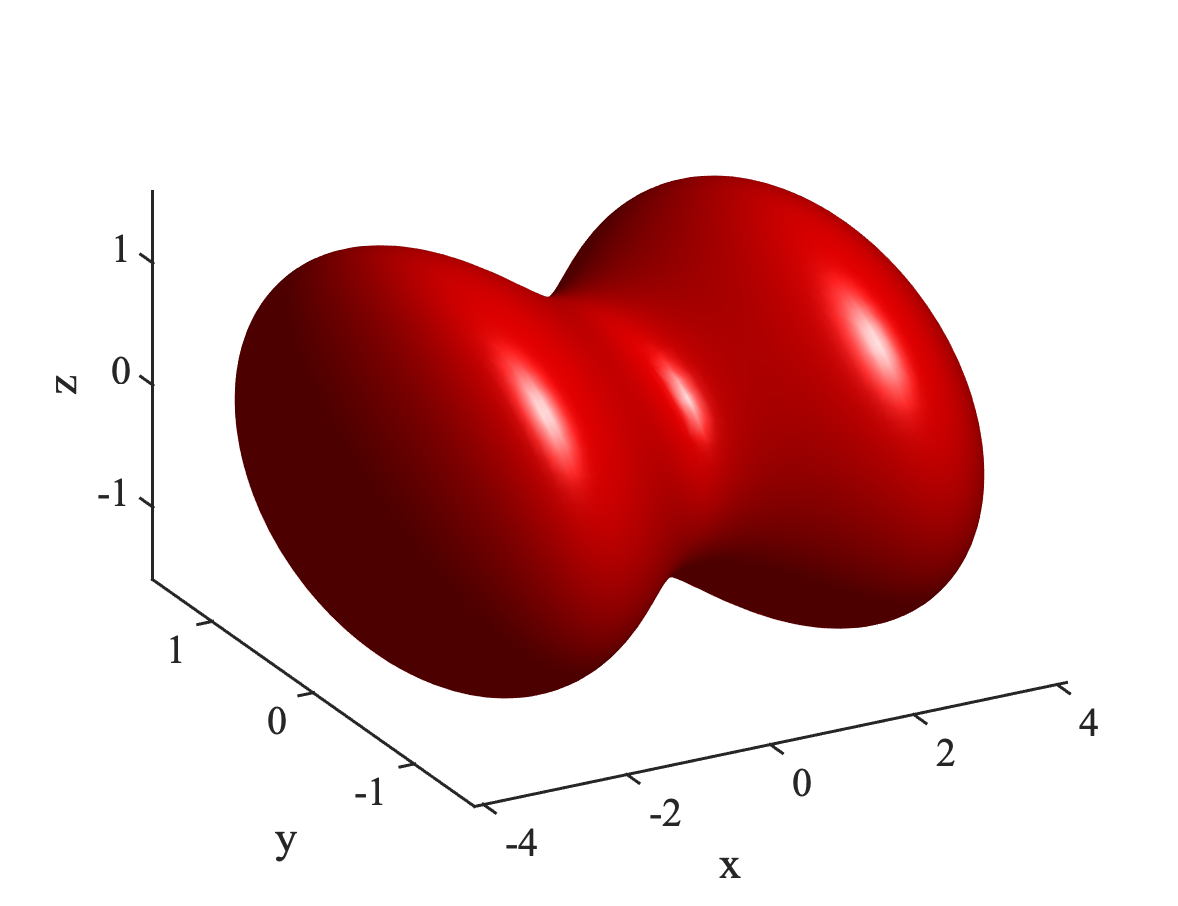}
\caption{The interface $\Gamma$.}
\label{fig: dumbbell}
\end{figure}

\begin{figure}[H]
\centering
\includegraphics[width=11cm,height=4cm]{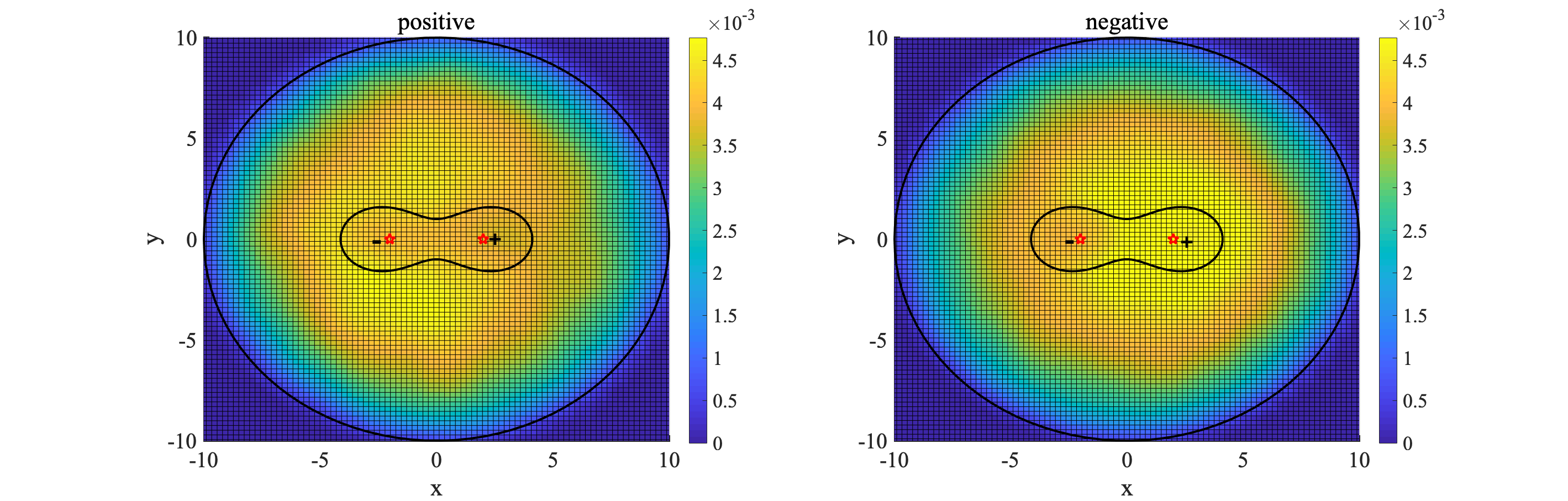}
\caption{The kernel density estimates in the $xOy$ planes. Left: positive charges; right: negative charges. The red pentagrams represent for the two singular free charges, the black curves denote the inner and outer boundaries of the truncated domain.}
\label{fig: xOy}
\end{figure}

\begin{figure}[H]
\centering
\includegraphics[width=11cm,height=4cm]{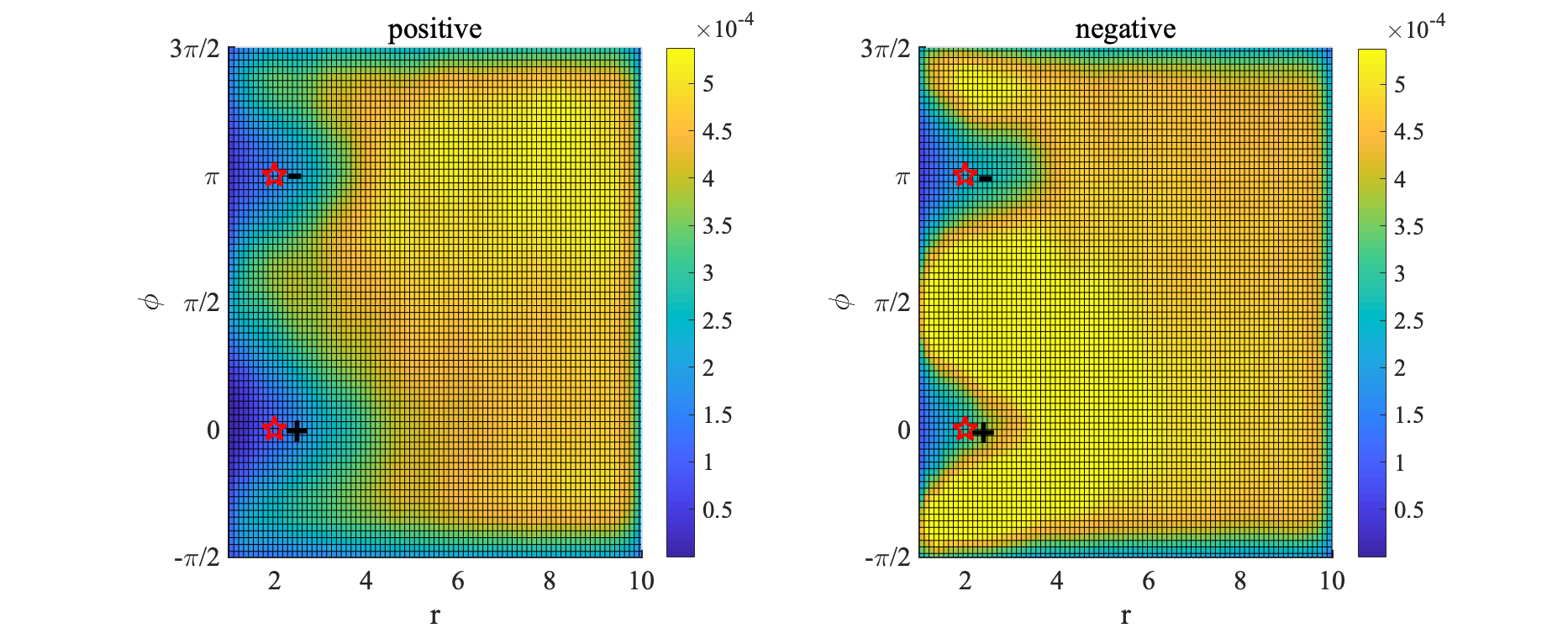}
\caption{The kernel density estimates in the $r-\phi$ plane. Left: positive charges; right: negative charges. The red pentagrams represent for the two singular free charges.}
\label{fig: r_phi}
\end{figure}

\begin{figure}[H]
\centering
\includegraphics[width=0.9\textwidth,height=0.4\textwidth]{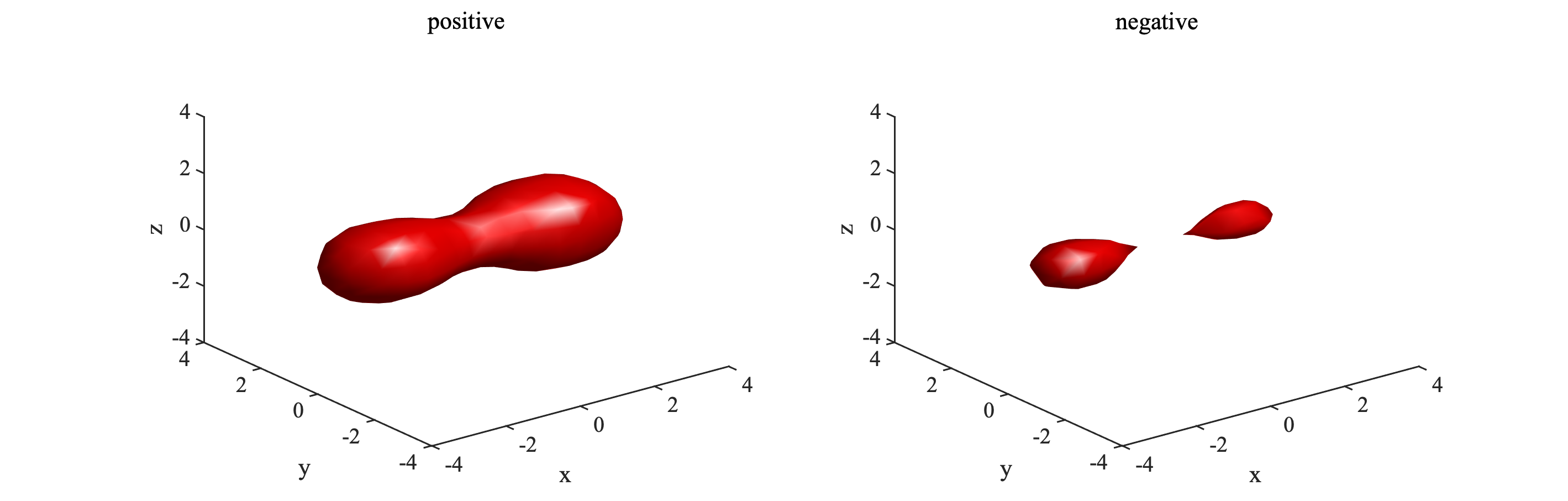}
\caption{The isosurface of density distribution. Left: $\rho_+=8e-4$; right: $\rho_-=8e-4$.}
\label{fig: isosurface}
\end{figure}

In this example, we apply time splitting RBM to the non-symmetric case. Consider $\mathcal{C}$ with a dumbbell-shaped interface
$$
\Gamma=\{(x,y,z)\in\mathbb{R}^3, (x^2+y^2+z^2)^2=a^2(x^2-y^2-z^2)+c\}.
$$
Take $a=4, c=17$, $\Gamma$ is shown in Figure \ref{fig: dumbbell}.
The free charge distribution is given by two singular charges
$$
\rho_f=\sum_{i=1}^2q_i\delta(\bm{x}-\bm{x}_i),\quad q_1=20, q_2=-1, \bm{x}_1=[2,0,0], \bm{x}_2=[-2,0,0].
$$
Again, we simulate particles inside truncated domain $\Omega_L=B_L\backslash \bar{\mathcal{C}}$ and collect samples from $100$ different time steps after the equilibrium to get more samples. Let $\nu=1, L=10, Q_+=10, q=10^{-3}, N_+=10^4, \tau=0.01$. Since $\Gamma$ is rotational symmetry about the $x$-axis, $\bm{x}_1, \bm{x}_2$ are located in the $x$-axis, then the positive and negative distributions are symmetric in the $yOz$ plane. The density distributions in the $xOy$ plane are in Figure \ref{fig: xOy} (similar to those in the $xOz$ plane). The same colorbar is used in the two subfigures. It is obvious that the density of positive charges is higher around $\bm{x}_2$ due to attraction, while lower around $\bm{x}_1$ due to repulsion. Meanwhile, the density of negative charges is higher around $\bm{x}_1$ and lower around $\bm{x}_2$. In the spherical coordinates $(r,\theta,\phi)$, one can observe non-symmetry in the $\phi$ direction, see Figure \ref{fig: r_phi}. As $\phi(\bm{x}_1)=0, \phi(\bm{x}_2)=\pi$, the negative charge density is lower around $\phi=\pi$ while higher around $\phi=0$. The attraction and repulsion are stronger around $\phi=0$ than $\phi=\pi$ because $|q_1|$ is much larger than $|q_2|$. Furthermore, in order to show the distribution intuitively, we plot the isosurface $\rho_{\pm}=8e-4$ in Figure \ref{fig: isosurface}. This example shows that random batch particle methods can be easily applied to complex geometry.

\subsection{Charge reversal in salty environment}
\label{subsec: colloid}
Though aiming at solving the PNP and PB equations using random batch particle simulations in this work, we emphasize that the random batch idea can be applied to other situations. In this example, we consider the case when a highly charged colloid is put into a solution containing trivalent counterions and monovalent coions. It has been observed in experiments \cite{experiment}, theories \cite{RevModPhys} and MD simulations \cite{colloid,RBE} that the effective charge of the colloid-ion complex becomes oppositely charged, which is the so-called charge reversal phenomenon. This is because the multivalent ions form a strongly correlated liquid on the surface of the colloid, so its charge can be overcompensated. Also, the inversion can not be predicted by the classical PB theory \cite{RevModPhys}.

Following the model setup in \cite{colloid}, we consider a spherical colloid of radius $R=50$ with a point charge $Q_f=-300$ at its center, in a solution of asymmetric $3:1$ salt, i.e., $z_{+}=3, z_{-}=-1$. There are $N_{+}=200$ trivalent counterions and $N_{-}=300$ monovalent coions, which are treated as uniform-sized with radius $R_o=4$. Thus, the total particle number $N=500$. The simulation domain is a spherical of radius $L=140$. The colloid is fixed in the center of the simulation domain. It is clear that the system is neutral, the colloidal surface charge density is $\sigma_f=\frac{1}{4\pi R^2}Q_f=0.0095$ due to symmetry. Particles of charges $q_i$ and $q_j$ interact via the Coulomb interaction $\phi_{ij}=\ell_b\frac{q_i q_j}{r_{ij}}$ with Bjerrum length $\ell_b=7.1$ (a number due to scaling using suitable units). Furthermore, the repulsive force is modeled via the shifted Lennard-Jones potential
\begin{equation*}
     \phi_{\mathrm{LJ}}=\left\{
     \begin{aligned}
  &4\epsilon\left[\left(\frac{\sigma}{r-r_{\text{off}}}\right)^{12}-\left(\frac{\sigma}{r-r_{\text{off}}}\right)^{6}\right]+
\epsilon,    \quad&  r-r_{\text{off}}<r_c,\\
      &0,  &\text{otherwise},
     \end{aligned}\right.
\end{equation*}
where the offset is $r_{\text{off}}=R+R_o$ between colloid and ion and $r_{\text{off}}=2R_o$ between ions. In the simulation, $\epsilon=1$, $\sigma=1$, $r_c=2^{1/6}\sigma$.

\begin{figure}[H]
\centering
\includegraphics[width=0.8\textwidth]{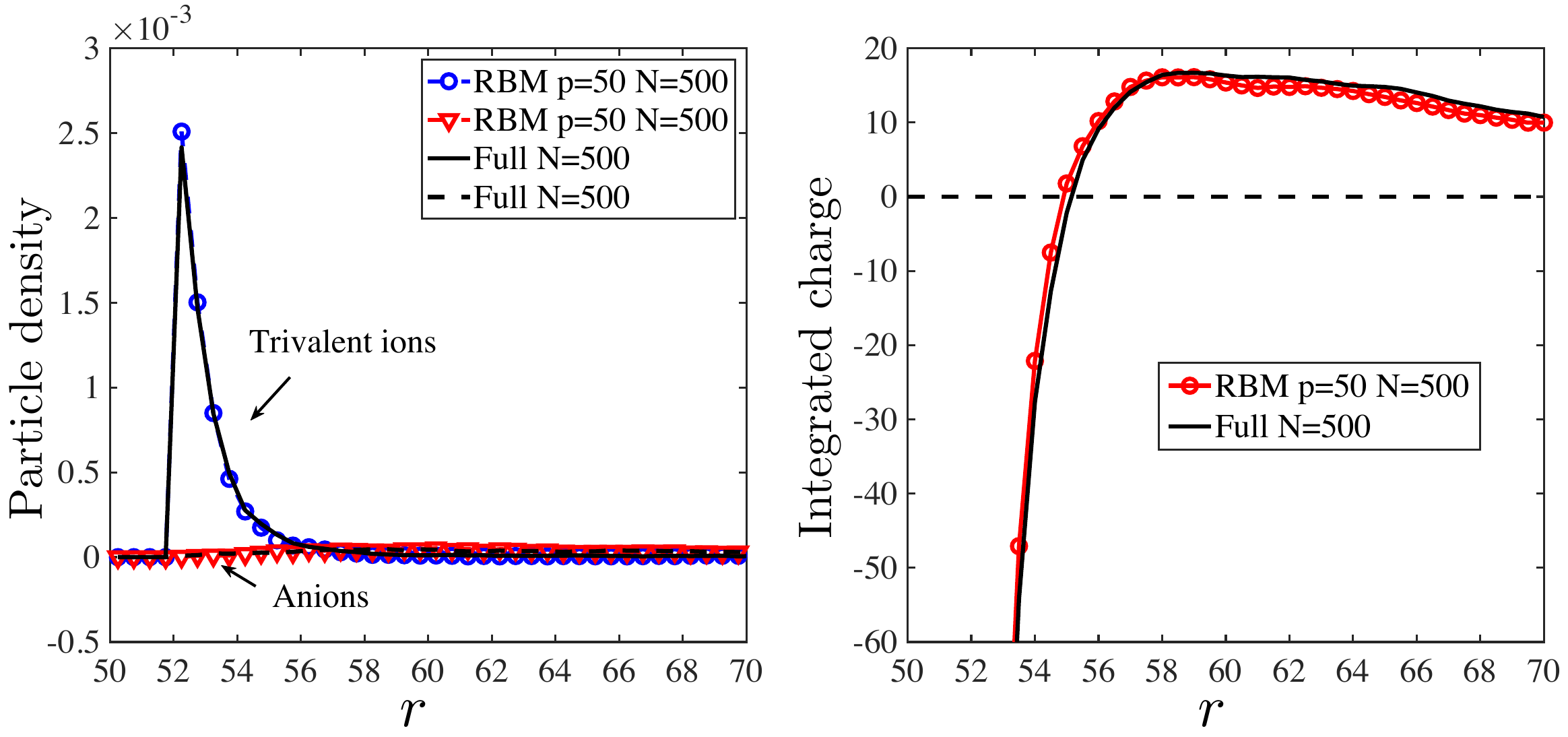}
\caption{Charge reversal test. Left: radial density distributions of different ion types; right: integrated charge against the distance $r$ from the colloid center.}
\label{fig: colloid}
\end{figure}

Let $\tau=0.005$, the results of random batch simulation with batch size $p=50$ are shown in Figure \ref{fig: colloid}, the full particle simulation is used as the reference solution. 
From the left subfigure, we can see that the distribution of trivalent counterions is highly peaked around the colloid due to strong attraction. In the right subfigure, the integrated charge distribution, which is the total charge within the distance $r$ from the colloidal center, is plotted. It clearly depicts that the charge of colloid is overcompensated up to $16$ in about two ion diameters from the colloidal surface. 
This tells that our random batch particle method can effectively capture the charge reversal phenomenon.

\begin{figure}[H]
\centering
\includegraphics[width=0.75\textwidth]{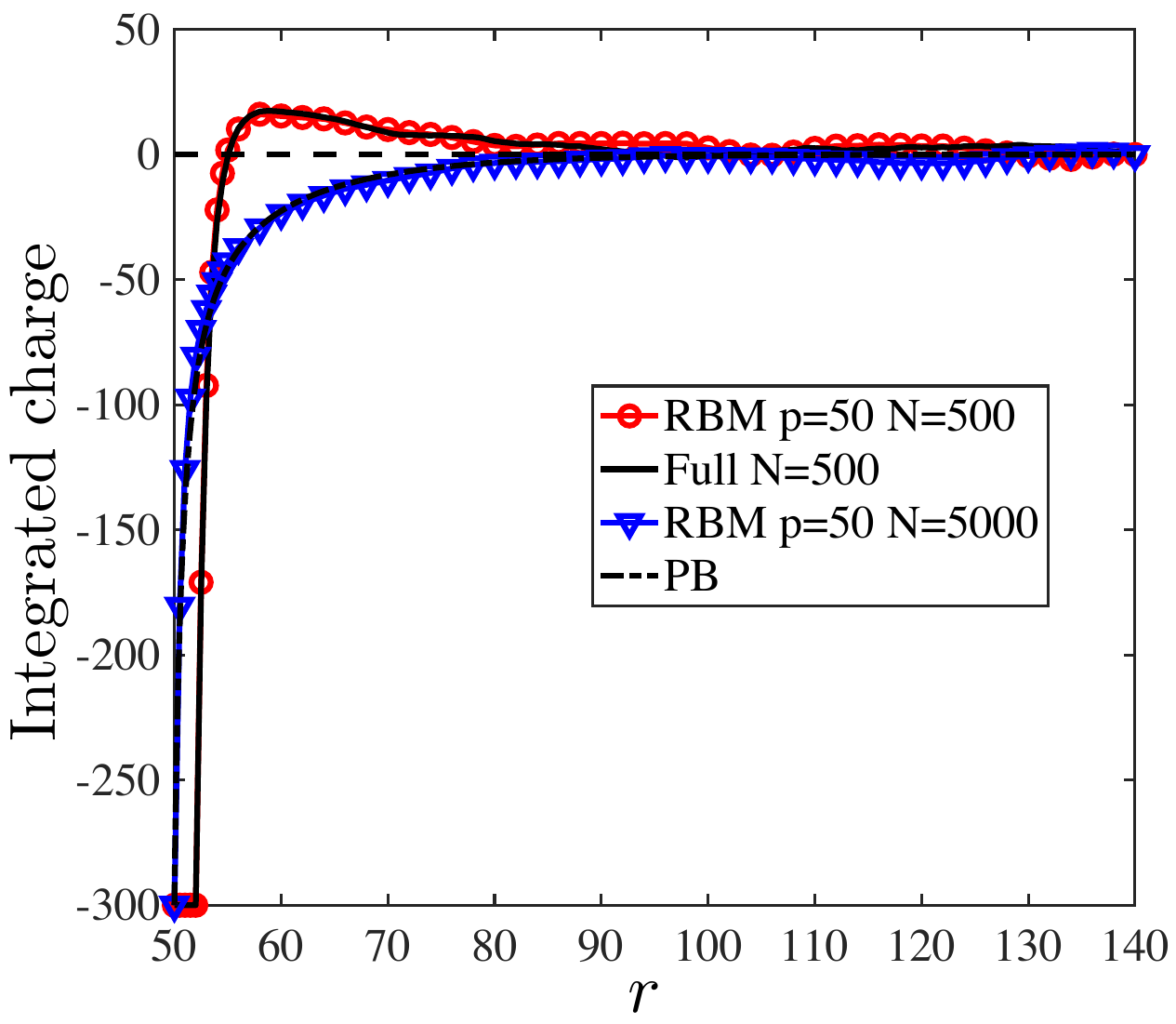}
\caption{Integrated charge against the distance $r$ from the colloid center. Comparison of RBM systems when $N=500$ and $N=5000$.}
\label{fig: transition}
\end{figure}

As a comparison, we increase the particle numbers $N_{\pm}$ so that the charge a particle carries is $q_{\pm}=Q_{\pm}/N_{\pm}\neq z_{\pm}$. In this case, the particle is numerical particle instead of physical particle. As $N_{\pm}\rightarrow \infty$, the ion size tends to zero, the system arrives at the mean-field regime, we can expect the simulation result gives an approximation to the PB equation of asymmetric $3:1$ salt. Namely, the following PB equation in the radial axis
\begin{equation}
    \label{eq: colloid}
     \left\{
    \begin{array}{lll}
    -\nu\left(\partial_{rr}\Phi+\frac{2}{r}\partial_r\Phi\right)=3\rho_{+}-\rho_{-}, \quad  r\in (R, L),\\[2mm]
    \rho_{+}=\rho_{\infty} e^{-3\Phi},\quad \rho_{-}=3 \rho_{\infty}e^{\Phi},\\[2mm]
    -\nu \partial_r\Phi(R)=\sigma_f,\quad \partial_r\Phi(L)=0.
    \end{array}\right.
\end{equation}
The parameter $\nu$ in \eqref{eq: colloid} is $1/(4\pi\ell_b)$ after rescaling. Increasing the particle number by $10$ times, we can observe a transition from the charge reversal phenomenon to the monotone PB regime in Figure \ref{fig: transition}. The curves for $N=500$ ($N_+=200, N_-=300$) are exactly those in the second plot of Figure \ref{fig: colloid}. When $N=5000$, the ions are assumed volumeless, i.e., $R_o=0$, the parameters of $\phi_{\mathrm{LJ}}$ are chosen as $ \epsilon=1/100, \sigma=1/\sqrt[3]{10}, r_c=2^{1/6}\sigma$ (we roughly took the lengths by scaling like $N^{-1/3}$ while $\epsilon$ by scaling like $q_{\pm}^2$ or $N^{-2}$), the FD solution to the PB equation \eqref{eq: colloid} is used as reference.


This example shows that our random batch particle methods inherit the following advantages of the MD simulations. First, some additional concerned physical effects like the Lennard-Jones interaction can be conveniently incorporated into our methods, which is not easy to be done in the PDE model. Second, while our methods can be used as numerical methods for the PNP and PB equations in the mean-field scaling, they can capture some interesting physics in MD simulations (like the charge reversal phenomenon) when the scaling is far away from the mean-field regime. Hence, our methods are far beyond being numerical methods for the PNP and PB equations.

\section{Conclusions}
\label{sec: Conclusions}

In this paper, we proposed some random batch particle methods with $\cO(N)$ cost per time step for the PNP and PB equations using RBM. 
By direct particle simulation, we can keep track of the dynamics of the particle as well as the equilibrium.
Meanwhile, the particle methods are insensitive to dimension and applicable to complex geometry. 
We do the particle simulation in a truncated domain and give an estimate saying that such a truncation makes sense for the PB equation in the symmetric case. Two splitting strategies are given to deal with the singular Coulomb potential. Also, it is feasible to take Lennard-Jones potential into account for physical purpose.

Though we mainly illustrated the methods using the $1:1$ electrolytes, the particle methods can be easily extended to solutions with multi-ionic species. Besides, it is interesting to do simulation for the problem with different dielectric constants in and outside the interface $\Gamma$, as discussed in Appendix \ref{sec: dielectric}.

\section*{Acknowledgement}
This work is partially supported by the National Key R\&D Program of China, Project Number 2021YFA1002800. The work of L. Li was partially sponsored by the Strategic Priority Research Program of Chinese Academy of Sciences, Grant No. XDA25010403, and NSFC 11901389, 12031013. The work of J.-G. Liu was supported by NSF DMS-2106988. Y. Tang would like to thank Prof. Shi Jin for his guidance and helpful discussions. The authors are very grateful to the anonymous referees for careful reviewing. The insightful suggestions have helped to improve our article.

\appendix
\section{Discussion about different dielectric constants}
\label{sec: dielectric}
In this section, we discuss the effects if the dielectric constant $\varepsilon$ inside $\mathcal{C}$ is different from the one outside.

First of all, recall the classical electrodynamics theory \cite{Griffiths2017}. By Gauss's law, 
$$
\varepsilon_0 \nabla \cdot E = \rho_{F}+\rho_p, \quad \rho_p=-\nabla \cdot P,
$$
where $\varepsilon_0$ is the dielectric constant in the vacuum, $P$ is the induced electric field, $E$ is the total electric field, and $\rho_F$, $\rho_p$ and $\rho_F+\rho_p$ represent for the free charge distribution, the induced charge distribution and the total charge distribution respectively.
Introducing the electric displacement vector
$$
D=\varepsilon_0E+P=\varepsilon E
$$
yields
$$
\nabla \cdot D=\nabla\cdot (\varepsilon E)=\rho_{F}.
$$
In other words, the total field $E$ is related to the free charge distribution if we put the varying electricity inside. The effects of induced charge are included in $\varepsilon$.

Suppose that $\varepsilon$ is piece-wise constant inside and outside $\mathcal{C}$, i.e.
 \begin{equation*}
 \varepsilon(\bm{x})=
 \left\{
 	\begin{array}{ll}
 		\varepsilon_1,\quad \bm{x}\in \mathcal{C},\\
 		\varepsilon_2,\quad \bm{x}\in \Omega.\\
 	\end{array}
 	\right.
 \end{equation*} The interface condition on $\Gamma$ is given by 
\begin{equation}
\label{eq: jump}
\bm{n}\times (E_2-E_1)=0,\quad \bm{n}\cdot (D_2-D_1)=0,
\end{equation}
where we assume there is no free surface charge on the interface.
Here, $E_1, E_2$ and $D_1, D_2$ are the total electric field and the electric displacement on the interior and exterior side of the cell respectively; $\bm{n}$ is the unit normal vector on the interface $\Gamma$ pointing to $\Omega$. This means the tangential components of $E$ are continuous on two sides of the matter, while the normal components of $E$ satisfy $\varepsilon_1 E_1 \cdot \bm{n} = \varepsilon_2 E_2\cdot \bm{n}$.
Hence, the Poisson equation in the entire space for our problem reads
\begin{equation}
\label{eq: Poisson entire}
\left\{
\begin{array}{ll}
	-\varepsilon_1\Delta \Phi=e\rho_f,\quad \bm{x} \in \mathcal{C}\\[2mm]
	-\varepsilon_2\Delta \Phi=e\sum\limits_{j=1}^J z_j\rho_j,\quad \bm{x} \in \Omega,\\[3mm]
	[\Phi]\big|_{\Gamma}=0,\quad [\varepsilon(\bm{x})\nabla\Phi\cdot \bm{n}]\big|_{\Gamma}=0.
	\end{array}
	\right.
\end{equation}
where $\rho_f$ is supported in $\mathcal{C}$ and $\rho_j$ is supported in $\Omega$.
In order to solve \eqref{eq: Poisson entire}, one needs to find the fundamental solution of
\begin{equation}\label{eq: fundamental sol of e(x)}
	-\nabla\cdot(\varepsilon(\bm{x})\nabla \Phi)=\delta (\bm{x}-\bm{x}_0),\quad \bm{x}_0\in \Omega.
\end{equation}
Clearly, the fundamental solution in $\Omega$ is not given by
\begin{equation*}
	\Phi_2(\bm{x})=\left\{
	\begin{array}{lll}
-\frac{1}{2\pi\varepsilon_2}\ln |\bm{x}-\bm{x}_0|,&\quad d=2,\\[3mm]
\frac{1}{d(d-2)\alpha(d)\varepsilon_2|\bm{x}-\bm{x}_0|^{d-2}},&\quad d\geq 3,
	\end{array}
    \right.
\end{equation*}
any more. Suppose otherwise.
It follows from $[\varepsilon(\bm{x})\nabla\Phi\cdot \bm{n}]\big|_{\Gamma}=0$ that
the solution to \eqref{eq: fundamental sol of e(x)} in $\mathcal{C} $ is 
\begin{equation*}
	\Phi_1(\bm{x})=\left\{
	\begin{array}{lll}
-\frac{1}{2\pi\varepsilon_1}\ln |\bm{x}-\bm{x}_0|+C,&\quad d=2,\\[3mm]
\frac{1}{d(d-2)\alpha(d)\varepsilon_1|\bm{x}-\bm{x}_0|^{d-2}}+C,&\quad d\geq 3,
	\end{array}
    \right.
\end{equation*}
where $C$ is an arbitrary constant. However, a simple computation shows that there does not exist a $C$ such that $[\Phi]=0$ on the interface $\Gamma$ if $\varepsilon_1\neq \varepsilon_2$.
Therefore, when we apply RBM for variable $\varepsilon$ cases directly, the interaction between two charges outside the cell no longer has a simple formula. 
In this case, we may use the image charge method to compute the interaction. For example, in \cite{Cai09}, many dielectric objects with different geometries are considered. In the case that the geometry is not a ball, the free charge $\rho_f$ in $\mathcal{C}$ (macromolecule or cell) can be modeled by multiple balls and the corresponding image charge method has been discussed in \cite{PhysRevE.87.013307}.

\section{Some missing proofs}
\label{sec: proofs}
\begin{proof}[Proof of Lemma \ref{lem: comparison principle}]
Given $\bm{x}\in \Omega$.
	Since $\lim\limits_{|\bm{x}|\rightarrow \infty} u(\bm{x})=0$, we have
	$$\forall \ \epsilon>0 \ \exists \ L \ \text{sufficiently large}, \quad s.t. \quad \bm{x}\in B_L,\quad |u(\bm{x}^{\prime})|<\epsilon\quad \text{for} \ |\bm{x}^{\prime}| \in \partial B_L.$$
	Let $\Omega_L=B_{L}\backslash \bar{\mathcal{C}}$. Then
	\begin{equation*}
	\left\{
\begin{array}{lll}
    \mathcal{L}u\leq 0,\quad \text{in} \ \Omega_L,\\[2mm]
	-\dfrac{\partial u}{\partial \bm{n}}\leq 0,\quad \text{on} \ \Gamma,\\[3mm]
	|u|<\epsilon, \quad \text{on} \ \partial B_{L}.
\end{array}
\right.
\end{equation*}
	By the maximum principle, $u$ attains its nonnegative maximum on $\partial \Omega_L$.

	If $u$ attains its nonnegative maximum on $\bm{x}_0 \in \Gamma$, 
	apply Hopf's lemma, $-\dfrac{\partial u}{\partial \bm{n}}\Big|_{\bm{x}=\bm{x}_0}> 0$, this contradicts with $-\dfrac{\partial u}{\partial \bm{n}}\Big|_{\Gamma}\leq 0$.
    Otherwise, $u$ attains its nonnegative maximum on $\partial B_{L}$, 
    so $u(\bm{x})\leq \left|u\right|_{\partial B_{L}}<\epsilon$.

    By the arbitrariness of $\epsilon$, $u\leq 0.$
\end{proof}

\begin{proof}[Proof of Proposition \ref{prop: decay rate} in 3D]
Throughout this proof, we will denote any generic constants as $C$, which may change from line to line.

STEP 1:
	Introduce a function of $s\in \mathbb{R}$,
	$$
	p(s)=\left\{\begin{array}{ll}
	\vspace{2mm}\dfrac{e^s-e^{-s}}{s}\rho_{\infty},&\quad s\neq 0,\\
	2\rho_{\infty},&\quad s=0.
	\end{array}\right.
	$$
	Then $p(s)\geq 2\rho_{\infty}>0$ is a continuous even function.

	So \eqref{eq:PB sym} can be rewritten as 
$$ \left\{
\begin{array}{ll} 
	-\nu \Delta \Phi+ p(\Phi)\Phi=0,\quad \bm{x}\in \Omega,\\[3mm]
	-\dfrac{\partial \Phi}{\partial \bm{n}}\Big|_{\Gamma} =\sigma_f,\quad \Phi(\bm{x}) \rightarrow 0 \ \text{as} \ |\bm{x}| \rightarrow \infty.
\end{array}
\right.
$$
Denote $\Sigma_f=\|\sigma_f\|_{L^{\infty}(\Gamma)}< +\infty$ and consider the linear problem 
\begin{equation}
\label{eq:linear PB}
	\left\{
\begin{array}{ll} 
	-\nu \Delta \bar{\Phi} + 2\rho_{\infty} \bar{\Phi}=0,\quad \bm{x}\in \Omega,\\[3mm]
	-\dfrac{\partial \bar{\Phi}}{\partial \bm{n}}\Big|_{\Gamma} =\Sigma_f,\quad \bar{\Phi}(\bm{x}) \rightarrow 0 \ \text{as} \ |\bm{x}| \rightarrow \infty.
\end{array}
\right.
\end{equation}
By Lemma \ref{lem: comparison principle}, $\bar{\Phi}\geq 0$.

Let $u^{\pm}=\bar{\Phi}\pm\Phi$, then 
$$ \left\{
\begin{array}{ll} 
	-\nu \Delta u^{\pm}+ p(\Phi)u^{\pm}=\left(p(\Phi)-2\rho_{\infty}\right)\bar{\Phi}\geq 0,\quad \bm{x}\in \Omega,\\[3mm]
	-\dfrac{\partial u^{\pm}}{\partial \bm{n}}\Big|_{\Gamma} =\Sigma_f\pm\sigma_f\geq 0,\quad u^{\pm} \rightarrow 0 \ \text{as} \ |\bm{x}| \rightarrow \infty.
\end{array}
\right.
$$
By Lemma \ref{lem: comparison principle}, $u^{\pm}\geq 0$. That is,
$$|\Phi|\leq \bar{\Phi}, \quad \bm{x}\in \Omega.$$

STEP 2: 
 In order to express the solution $\bar{\Phi}$ of the exterior Neumann problem \eqref{eq:linear PB}, first construct the interior Dirichlet problem
$$
\left\{
\begin{array}{ll}
\vspace{2mm}-\nu \Delta \bar{\Phi} + 2\rho_{\infty}\bar{\Phi}=0 \quad \text{in}\ \mathcal{C},\\[2mm]
\bar{\Phi}=\bar{\Phi}^e \quad \text{on} \ \Gamma.
\end{array}\right.
$$
Recall that for $\bm{x}\in \Gamma$, $\bar{\Phi}^e(\bm{x})$ and $\bar{\Phi}^i(\bm{x})$ represent the limits from exterior and interior side.
Hence, $[\bar{\Phi}]=0$. Then, by the representation formula \eqref{eq: representation formula}, we have
\begin{equation}
\label{eq: phi1}
\bar{\Phi}(\bm{x})=\int_{\Gamma}q(\bm{y})G(\bm{x}-\bm{y})\mathrm{d} S_{\bm{y}},
\end{equation}
where $q=\left[\dfrac{\partial \bar{\Phi}}{\partial\bm{n}}\right]$ and
\begin{equation}\label{eq:Green}
G(\bm{x})=\dfrac{e^{-\sqrt{2\rho_{\infty}/\nu}|\bm{x}|}}{4\pi|\bm{x}|}
\end{equation}
is the three-dimensional Green function.

As for $\bm{x}\in \Gamma,\quad \dfrac{\partial \bar{\Phi}^e}{\partial \bm{n}}(\bm{x})=\dfrac{\frac{\partial \bar{\Phi}^e}{\partial \bm{n}}(\bm{x})+\frac{\partial \bar{\Phi}^i}{\partial \bm{n}}(\bm{x})}{2}+\dfrac{\frac{\partial \bar{\Phi}^e}{\partial \bm{n}}(\bm{x})-\frac{\partial \bar{\Phi}^i}{\partial \bm{n}}(\bm{x})}{2},$
$q$ satisfies
$$
-\Sigma_f=\int_{\Gamma}q(\bm{y})\dfrac{\partial G(\bm{x}-\bm{y})}{\partial \bm{n}_{\bm{x}}}\mathrm{d}S_{\bm{x}}-\frac{q(\bm{x})}{2}.
$$
Define a kernel
$k: \Gamma \times \Gamma \rightarrow \mathbb{R}$ as
\begin{equation}
\label{eq:kernel}
	k(\bm{x},\bm{y})=2\dfrac{\partial G(\bm{x}-\bm{y})}{\partial \bm{n}_{\bm{x}}}
\end{equation}
and define an integral operator $K$ by
\begin{equation}
\label{eq: operator} 
	Kf(\bm{x})=\int_{\Gamma} k(\bm{x},\bm{y})f(\bm{y})\mathrm{d}S_{\bm{y}},\quad \bm{x}\in\Gamma .
\end{equation}
Then we are led to the Fredholm equation of the second kind:
\begin{equation}
	\label{eq: Fredholm}
	(I-K)q=2\Sigma_f.
\end{equation}

From Lemma \ref{lem: compactness of K}, we know $K$ is compact as operators in $ L^2(\Gamma)\rightarrow L^2(\Gamma)$ and also in $C(\Gamma)\rightarrow C(\Gamma)$. $I-K$ is a Fredholm operator. Therefore, given $\Sigma_f \in L^2(\Gamma)$, there exists a unique $q\in L^2(\Gamma)$ such that $(I-K)q=2\Sigma_f$. Besides, since $K$ is compact in $ C(\Gamma)\rightarrow C(\Gamma)$, $q\in C(\Gamma)$. Hence, $\bar{\Phi}\in C^2(\Omega)\bigcap C(\bar{\Omega})$ from \eqref{eq: phi1}.

Now we are in a position to show the exponentially decay of $\bar{\Phi}$.
For all $\bm{x} \in \Omega$ satisfying $|\bm{x}|\geq 2\sup\limits_{\bm{y}\in\Gamma} |\bm{y}|$, it follows from \eqref{eq: phi1} that
$$
\left|\bar{\Phi}(\bm{x})\right|\leq\|q\|_{L^2(\Gamma)}\left(\int_{\Gamma}\frac{e^{-2\sqrt{2\rho_{\infty}/\nu}}|\bm{x}-\bm{y}|}{4^2\pi^2|\bm{x}-\bm{y}|^2}\mathrm{d}S_{\bm{y}}\right)^{\half}
\leq\|q\|_{L^2(\Gamma)}\sqrt{|\Gamma|}\dfrac{e^{-\sqrt{\rho_{\infty}/2\nu}|\bm{x}|}}{2\pi|\bm{x}|},
$$
where $|\Gamma|$ is the area of $\Gamma$.
Thus, there exists $ R_1=2\sup\limits_{\bm{y}\in\Gamma} |\bm{y}|$ such that
$$\left|\Phi(\bm{x})\right|\leq \left|\bar{\Phi}(\bm{x})\right| \leq \frac{C}{|\bm{x}|}e^{-\frac{C}{\sqrt{\nu}}|\bm{x}|},\quad |\bm{x}|>R_1.$$
Besides, $\Phi\in L^{\infty}(\Omega)$ as $\bar{\Phi}\in C(\bar{\Omega})$. 

STEP 3: 
Look back to the original PB equation \eqref{eq:PB sym}
$$-\nu \Delta \Phi + 2\rho_{\infty} \Phi=2\rho_{\infty} (\Phi-\sinh \Phi),\quad \bm{x}\in \Omega.$$
Let $f:=2\rho_{\infty} (\Phi-\sinh \Phi)\in L^{\infty}(\Omega)$. Then, since $\Phi \in L^{\infty}(\Omega)$, there exists $R_2$ large enough such that
\begin{equation}\label{eq:est f}
\left|f(\bm{x})\right| \leq C\Phi^2\leq \frac{C}{|\bm{x}|^2}e^{-\frac{C}{\sqrt{\nu}}|\bm{x}|},\quad |\bm{x}|>R_2.
\end{equation}

Next, by similar arguments with STEP 2, we show the exponentially decay property of $\nabla \Phi$.
Construct the interior Dirichlet problem
$$
\left\{
\begin{array}{ll}
\vspace{2mm}-\nu \Delta \Phi + 2\rho_{\infty}\Phi=0 \quad \text{in}\ \mathcal{C},\\[2mm]
\Phi=\Phi^e \quad \text{on} \ \Gamma.
\end{array}\right.
$$
Then $[\Phi]=0$ and by the representation formula \eqref{eq: representation formula}, we have
\begin{equation}
\label{eq: phi2}
\Phi(\bm{x})=\int_{\Gamma}\tilde{q}(\bm{y})G(\bm{x}-\bm{y})\mathrm{d}S_{\bm{y}}+\frac{1}{\nu}\int_{\Omega}G(\bm{x}-\bm{y})f(\bm{y})\mathrm{d}\bm{y},
\end{equation}
where $\tilde{q}=\left[\dfrac{\partial \Phi}{\partial\bm{n}}\right]$.
Similarly, we are led to the Fredholm equation of the second kind
\begin{equation}
\label{eq:Fredholm2}
	(I-K)\tilde{q}=2(g+\sigma_f)
\end{equation}
with $g=\frac{1}{\nu}\int_{\Omega}\frac{\partial G(\bm{x}-\bm{y})}{\partial \bm{n}_{\bm{x}}}f(\bm{y})\mathrm{d}\bm{y}\in L^2(\Gamma)$.
Therefore, there exists a unique $\tilde{q}\in L^2(\Gamma)$ satisfing \eqref{eq:Fredholm2}.

Now we estimate the gradient of $\Phi$. One can deduce from \eqref{eq: phi2} that
\begin{align}
\nabla\Phi(\bm{x})=&\frac{1}{\nu}\int_{\Omega}\nabla_{\bm{x}}G(\bm{x}-\bm{y})f(\bm{y})\mathrm{d}\bm{y}+\int_{\Gamma}\nabla_{\bm{x}}G(\bm{x}-\bm{y})\tilde{q}(\bm{y})\mathrm{d}S_{\bm{y}}.
\label{eq: phi3}
\end{align}
According to \eqref{eq:Green}, one has
\begin{equation}\label{eq:grad G}
\left|\nabla_{\bm{x}} G(\bm{x}-\bm{y})\right|\leq \frac{\sqrt{2\rho_{\infty}/\nu}|\bm{x}-\bm{y}|+1}{4\pi|\bm{x}-\bm{y}|^2}e^{-\sqrt{2\rho_{\infty}/\nu}|\bm{x}-\bm{y}|}.
\end{equation}
For $\bm{x} \in \Omega$ satisfying $|\bm{x}|\geq 2R_2$, consider 
\begin{align*}
\uppercase\expandafter{\romannumeral1}:=&\int_{\Omega}\big|\nabla_{\bm{x}} G(\bm{x}-\bm{y})f(\bm{y})\big|\mathrm{d}\bm{y}\\[2mm]
=&\left(\int_{\Omega_1}+\int_{\Omega_2}+\int_{\Omega_3}\right)\big|\nabla_{\bm{x}} G(\bm{x}-\bm{y})f(\bm{y})\big|\mathrm{d}\bm{y}:=I_1+I_2+I_3,
\end{align*}
where
$\Omega_1=\{ |\bm{x}|/2 \leq |\bm{y}|\leq 2|\bm{x}|\}, \Omega_2=\{ |\bm{y}|\geq 2|\bm{x}|\}, \Omega_3=\{|\bm{y}|\leq |\bm{x}|/2\}\cap \Omega$. Then, by \eqref{eq:est f} and \eqref{eq:grad G}, one can obtain that
\begin{align*}
\uppercase\expandafter{\romannumeral1}_1
\leq& \|f\|_{L^{\infty}(\Omega_1)}\int_{\Omega_1}\left|\nabla_{\bm{x}} G(\bm{x}-\bm{y})\right|\mathrm{d}\bm{y} \leq 
\frac{C}{|\bm{x}|^2}e^{-\frac{C}{\sqrt{\nu}}|\bm{x}|}\left(2\sqrt{\nu}-(2\sqrt{\nu}+C|\bm{x}|)e^{-\frac{C}{\sqrt{\nu}}|\bm{x}|}\right)\leq \frac{C\sqrt{\nu}}{|\bm{x}|^2}e^{-\frac{C}{\sqrt{\nu}}|\bm{x}|},\\[2mm]
\uppercase\expandafter{\romannumeral1}_2
\leq& \|f\|_{L^1(\Omega_2)}\sup_{\Omega_2}\left|\nabla_{\bm{x}} G(\bm{x}-\bm{y})\right|
\leq C\sqrt{\nu}e^{-\frac{C}{\sqrt{\nu}}|\bm{x}|}\frac{\frac{C}{\sqrt{\nu}}|\bm{x}|+1}{|\bm{x}|^2}e^{-\frac{C}{\sqrt{\nu}}|\bm{x}|}\leq \frac{C\sqrt{\nu}}{|\bm{x}|^2}e^{-\frac{C}{\sqrt{\nu}}|\bm{x}|},\\[2mm]
\uppercase\expandafter{\romannumeral1}_3
\leq& \|f\|_{L^1(\Omega_3)}\sup_{\Omega_3}\left|\nabla_{\bm{x}} G(\bm{x}-\bm{y})\right|\leq C\sqrt{\nu}\left(1-e^{-\frac{C}{\sqrt{\nu}}|\bm{x}|}\right)\frac{\frac{C}{\sqrt{\nu}}|\bm{x}|+1}{|\bm{x}|^2}e^{-\frac{C}{\sqrt{\nu}}|\bm{x}|}\leq \frac{C\sqrt{\nu}}{|\bm{x}|^2}e^{-\frac{C}{\sqrt{\nu}}|\bm{x}|}.
\end{align*}
For $ \bm{x}\in \Omega$ satisfying $|\bm{x}|\geq R_1$, it follows from \eqref{eq:grad G} that
$$J:=
\int_{\Gamma}\left|\nabla_{\bm{x}} G(\bm{x}-\bm{y})\tilde{q}(\bm{y})\right|\mathrm{d} S_{\bm{y}}\leq \|\tilde{q}\|_{L^2(\Gamma)}\left(\int_{\Gamma}|\nabla_{\bm{x}} G(\bm{x}-\bm{y})|^2\mathrm{d}S_{\bm{y}}\right)^{\half}\leq C\frac{\frac{C}{\sqrt{\nu}}|\bm{x}|+1}{|\bm{x}|^2}e^{-\frac{C}{\sqrt{\nu}}|\bm{x}|}.
$$ 
Therefore, there exists a large $R=\max\{R_1,2R_2\}$ such that for all $|\bm{x}|\geq R$, 
$$
\left|\nabla\Phi(\bm{x})\right|\leq \frac{I}{\nu}+J\leq \frac{C}{\sqrt{\nu}|\bm{x}|^2}e^{-\frac{C}{\sqrt{\nu}}|\bm{x}|}.
$$
\end{proof}

\begin{lemma}
\label{lem: compactness of K}
	For smooth boundary $\Gamma\subset \mathbb{R}^2$, the operator $K$ defined by \eqref{eq:Green}, \eqref{eq:kernel}, \eqref{eq: operator} is compact both in $B(L^2(\Gamma))$ and $B(C(\Gamma))$.
\end{lemma}
\begin{proof}
For compactness in $C(\Gamma)\rightarrow C(\Gamma)$, it is similar to the Laplacian case, and see \cite{Hackbusch1995Integral} for the latter.

Now consider the compactness in $ L^2(\Gamma)\rightarrow L^2(\Gamma)$:
$$k(\bm{x},\bm{y})=2\bm{n}(\bm{x})\cdot \nabla G(\bm{x}-\bm{y})=\frac{-\bm{n}(\bm{x})\cdot (\bm{x}-\bm{y})}{2\pi}\dfrac{\sqrt{2\rho_{\infty}/\nu}|\bm{x}-\bm{y}|+1}{|\bm{x}-\bm{y}|^3}e^{-\sqrt{2\rho_{\infty}/\nu}|\bm{x}-\bm{y}|}.$$
Since $\Gamma$ is smooth, $\bm{n}(\bm{x})$ is almost orthogonal to $(\bm{x}-\bm{y})$ for $\bm{y}$ close to $\bm{x}$, there exists a bounded continuous function $c(\bm{x},\bm{y})$ such that
$$-\bm{n}(\bm{x})\cdot (\bm{x}-\bm{y})=|\bm{x}-\bm{y}|^2c(\bm{x},\bm{y}).$$
Thus the kernel is weakly singular (behaves like $|\bm{x}-\bm{y}|^{-1}$ as $\bm{x} \rightarrow \bm{y}$).

For $\varepsilon>0$, let
$$\quad k_{\varepsilon}(\bm{x},\bm{y})=\frac{ c(\bm{x},\bm{y})}{2\pi}\frac{\sqrt{2\rho_{\infty}/\nu}|\bm{x}-\bm{y}|+1}{|\bm{x}-\bm{y}|+\varepsilon}e^{-\sqrt{2\rho_{\infty}/\nu}|\bm{x}-\bm{y}|},\quad \bm{x},\bm{y} \in\Gamma.$$
$$K_{\varepsilon}f(\bm{x})=\int_{\Gamma} k_{\varepsilon}(\bm{x},\bm{y})f(\bm{y})\mathrm{d} S_{\bm{y}},\quad f\in L^2(\Gamma).$$
Then $k_{\varepsilon}\in L^2\left(\Gamma\times \Gamma\right)$, so  $K_\varepsilon: L^2(\Gamma)\rightarrow L^2(\Gamma)$ is a Hilbert-Schmidt operator with norm $\|K_{\varepsilon}\|_{HS}=\|k_{\varepsilon}\|_{L^2}<\infty$, and hence compact.

Next, we estimate $\|K_{\varepsilon}-K\|_{B(L^2(\Gamma))}$. For smooth surface $\Gamma$, there exist $n$ local charts $D_i\subset \mathbb{R}^2$ and mapping $\varphi_i$, $i=1, \cdots, n$ such that  
$$\varphi_i\circ k_\varepsilon\sim \frac{C_{i}|u-v|+C_{i}'}{|u-v|+\varepsilon}=:\bar{k}_{i,\varepsilon},\quad \varphi_i\circ k\sim \frac{C_{i}|u-v|+C_{i}'}{|u-v|}=:\bar{k}_{i},\quad u, v\in D_i.$$ 
Since $f\in L^2(\Gamma)$,$ \bar{f}_{i}=\varphi_i\circ f\in L^2(D_i)$. Define
$$\bar{K}_{i,\varepsilon}\bar{f}_i(v)=\int_{D_i}\bar{k}_{i,\varepsilon}(u,v)\bar{f}_i(u)\mathrm{d}u,\quad \bar{K}_i\bar{f}_i(v)=\int_{D_i}\bar{k}_{i}(u,v)\bar{f}_i(u)\mathrm{d}u.$$
Let
 $$
 g_{i,\varepsilon}(u)=\frac{C_{i}|u|+C_{i}'}{|u|+\varepsilon}1_{B_i}(u),\quad
 g_{i}(u)=\frac{C_{i}|u|+C_{i}'}{|u|}1_{B_i}(u),
 $$
 with bounded sets $B_i\subset \mathbb{R}^2$ s.t. $D_i\subset B_i,~i=1,\cdots, n$. Then, 
 $$\bar{K}_{i,\varepsilon}\bar{f}_i=g_{i,\varepsilon}*1_{D_i}\bar{f}_i,\quad \bar{K}_{i}\bar{f}_i=g_{i}*1_{D_i}\bar{f}_i.$$
Using local charts, it suffices to consider the operators $\bar{K}_{i,\varepsilon}, \bar{K}_i$.
Applying Young's convolution inequality, one has
$$
	\|(\bar{K}_{i,\varepsilon}-\bar{K}_i)\bar{f}_i\|_{L^2(D_i)}\leq\|(g_{i,\varepsilon}-g_i)*1_{D_i}\bar{f}_i\|_{L^2(\mathbb{R}^2)}\leq\|g_{i,\varepsilon}-g_i\|_{L^1(\mathbb{R}^2)}\|\bar{f}_i\|_{L^2(D_i)}.
$$
This implies
$\|\bar{K}_{i,\varepsilon}-\bar{K}_i\|_{B(L^2(D_i))}\leq\|g_{i,\varepsilon}-g_i\|_{L^1(\mathbb{R}^2)}\rightarrow 0 $ as $\varepsilon\rightarrow 0$ by dominant convergence theorem. Therefore, 
$$\|K_{\varepsilon}-K\|_{B(L^2(\Gamma))}\leq \sum_{i=1}^n\|\bar{K}_{i,\varepsilon}-\bar{K}_i\|_{B(L^2(D_i))}
\rightarrow 0 \ \text{as} \ \varepsilon\rightarrow 0,$$
i.e. 
$K_{\varepsilon}\rightarrow K $ in $ B(L^2(\Gamma))$ as $\varepsilon\rightarrow 0$. Hence, $K: L^2(\Gamma)\rightarrow L^2(\Gamma)$ is a compact operator.
\end{proof}
\bibliographystyle{plain}
\bibliography{rsRef.bib}

\end{document}